\documentclass{article}

\usepackage{arxiv}

\usepackage[utf8]{inputenc} 
\usepackage[T1]{fontenc}    
\usepackage{hyperref}       
\usepackage{url}            
\usepackage{booktabs}       
\usepackage{amsfonts}       
\usepackage{amsthm,amsmath,amssymb}
\usepackage{mathrsfs}
\usepackage{nicefrac}       
\usepackage{microtype}      
\usepackage{lipsum}
\usepackage{graphicx}
\usepackage{caption}
\usepackage{subcaption}
\usepackage{algpseudocode}
\graphicspath{ {./images/} }

\title{Unified interface flux evaluation in a general discontinuous Galerkin spectral element framework}

\author{
 Boyang Xia \\
  Department of Engineering \\
  King's College London \\
  Strand, London, WC2R 2LS \\
  \texttt{boyang.xia@kcl.ac.uk} \\
   \And
 David Moxey \\
  Department of Engineering \\
  King's College London \\
  Strand, London, WC2R 2LS \\
  \texttt{david.moxey@kcl.ac.uk} \\
}

\begin{document}
\maketitle
\begin{abstract}
  High-order discontinuous Galerkin spectral element methods (DGSEM) have
  received growing attention and development in recent years, especially in the
  regime of computational fluid dynamics. The inherent flexibility of the DG
  approach in handling non-conforming interfaces, such as those encountered in
  moving geometries or $hp$-refinement, presents a significant advantage for
  real-world simulations. Despite the well-established mathematical framework of
  DG methods, practical implementation challenges persist to boost performance
  and capability. Most previous studies only focus on certain choices of element
  shape or basis type in a structured mesh, although they have demonstrated the
  capability of DGSEM in complex flow simulations. This work discusses low-cost
  and unified interface flux evaluation approaches for general spectral elements
  in unstructured meshes, alongside their implementations in the open-source
  spectral element framework, Nektar++. The initial motivation arises from the
  discretisation of Helmholtz equations by the symmetric interior penalty
  method, in which the system matrix can easily become non-symmetric if the flux
  is not properly evaluated on non-conforming interfaces. We focus on the
  polynomial non-conforming case in this work, but extending to the geometric
  non-conforming case is theoretically possible. Comparisons of different
  approaches, trade-offs, and performance of our initial matrix-free
  implementation are also included, contributing to the broader discourse on
  high-performance spectral element method implementations.
\end{abstract}

\keywords{ spectral/\textit{hp} element method, discontinuous Galerkin methods, symmetric interior penalty method, matrix-free methods, SIMD vectorisation}

\section{Introduction}
In computational fluid dynamics (CFD), high-order methods have gained popularity
in research and industrial applications in recent
years~\cite{wang_high-order_2013}. High-order methods, together with modern
high-performance computing facilities, play an important role in achieving
high-fidelity, multiscale resolved turbulent flow simulations in industrial
applications. From a numerical perspective, diffusion and dispersion are greatly
reduced at higher orders, meaning that such methods are ideally placed for
tracking energetic flow structures such as jetting
vortices~\cite{lombard-2016,buscariolo-2022} or modelling problems involving
high separation~\cite{slaughter-2023, mengaldo-2020} across long time- and
length-scales. The use of a finite element discretisation also means that they
retain the geometric accuracy required to model complex geometries and provide
localised refinement, as needed for e.g. wall-resolved large-eddy
simulations~\cite{moxey-2015a}. From a computing perspective, it is also
well-documented that these high-order methods are well-suited to modern CPUs and
GPUs due to higher operational intensity. For a desired level of solution error,
we can reduce the total runtime by adjusting the polynomial order and mesh size
accordingly. This has been highlighted in a range of popular frameworks,
including deal.II~\cite{kronbichler-2018}, Dune~\cite{bastian-2010},
PyFR~\cite{witherden-2014a,witherden2025pyfr} and MFEM~\cite{andrej2024high}, or
efforts such as the CEED project~\cite{kolev2021efficient}.

Among the variety of high-order finite element methods, the discontinuous
Galerkin (DG) method in particular is an area of rapid growth and
development. DG is similar to classic finite volume methods (FVM), which allow
discontinuity across element boundaries via a numerical flux and are broadly
applicable across a range of hyperbolic problems. However, DG discretisations
typically have a more compact stencil: that is, elements only need to
communicate with their direct adjacent elements, rather than multiple layers of
adjacencies, regardless of discretisation order~\cite{bassi1997high}. In this
manner, the treatment of complex non-conformal mesh interfaces in cases such as
local $hp$-refinement and moving geometries is unified and independent of
element order. Although the theoretical analysis of DG has been well studied in
the past few decades~\cite{arnold_unified_2002,cockburn_discontinuous_2012},
their practical implementation in code varies according to different choices of
the element type, basis, and target applications, requiring additional analysis
and design. We are especially interested in these details, as they shape code
design as well as attainable performance.

Although the theoretical analysis of DG has been well studied in
the past few decades~\cite{arnold_unified_2002,cockburn_discontinuous_2012},
their practical implementation in code varies according to different choices of
the element type, basis, and target applications, requiring additional analysis
and design. We are especially interested in these details, as they shape code
design as well as attainable performance.

The most popular element shapes in spectral element methods are quadrilateral
and hexahedral elements, also known as tensorial elements since their expansions
and quadratures are constructed via tensor product. Such elements possess
several advantages. For instance, tensor computing techniques such
as sum-factorisation~\cite{orszag_spectral_1980} and collocation can be applied 
to reduce evaluation cost and enable matrix-free operators. 
There have been many encouraging developments towards
practical, robust, and efficient fluid dynamics solvers based on hexahedral
elements in recent years~\cite{fehn_stability_2017, fehn_robust_2018,
  bastian_matrix-free_2019, fehn_hybrid_2020, ferrer_horses3d_2023,
  jansson_neko_2024}, enabling larger-scale and more complex flow simulations.

The motivation of this work arises from the fact that the generation of
high-quality, fully hexahedral meshes suitable for the simulation of realistic
industrial geometries is an open and complex problem~\cite{pietroni2022hex}. In
these cases, the most common mesh decompositions involve either fully
tetrahedral or mixed prismatic-tetrahedral meshes. For tetrahedra, prisms and
pyramids, general tensorial expansions such as those presented in
in~\cite{karniadakis_spectralhp_1999} and~\cite{warburton_galerkin_1999} are
still possible. Although they result in higher computing cost relative to
hexahedra, we can still benefit from matrix-free implementations and achieve
over 50\% of the peak performance of modern CPUs, by making good use of
multi-level caches, vectorisation, and sum-factorisation in continuous Galerkin
methods~\cite{moxey_efficient_2020}..

A secondary interest is non-conforming interfaces, which appear when adjacent
elements are either misaligned (geometrical non-conforming) or have different
orders (polynomial non-conforming). The most well-known approach to handle such
cases is called the mortar element method~\cite{kopriva_staggered-grid_1998},
first considered for $C^0$-continuous spectral element
methods~\cite{maday,mavriplis}, but can also be applied to
DG~~\cite{kopriva_computation_2002,
  laughton_comparison_2021,heinz_high-order_2023,
  durrwachter_efficient_2021}. $L^2$ projection between local elements and
mortar elements is essential to retain convergence. In all previous
publications, the projection is performed between the solution degrees of
freedom, so that the inverse mass matrix is explicitly presented in the
formulation, which is typically expensive to compute. Constructing mortar
elements is cumbersome between unstructured triangular mesh interfaces; therefore, 
in these cases, we prefer point-to-point interpolation~\cite{laughton_comparison_2021}, 
which directly evaluates the adjacent solution at local points and 
computes the flux term locally. However, it is known to be sub-optimal in accuracy.

Our goal in this work is to address these gaps, by presenting a complete picture
of formulation, implementation design and performance analysis for the SIPG
discretisation of a Helmholtz operator on various element types. The paper is
structured as follows. In section~\ref{sec:basics}, we introduce the formulation
of the method and its existing implementation within Nektar++. We then introduce
an optimised workflow to evaluate any DG formulation, followed by general
interface flux evaluation designs. Section~\ref{sec:SIPG} presents a theoretical
analysis of the interior penalty method from a unique perspective, to illustrate
how flux evaluation on non-conforming interfaces can affect the symmetry of
systems. Based on the analysis, we propose two ways to handle non-conforming
interfaces in section~\ref{sec:unified}: one is to use sufficient quadrature,
and the other is to define a shared trace space and unify the quadrature for
both sides, similar to mortar elements. We discuss in detail how these two
approaches can be implemented to achieve unified, low-cost flux evaluations for
general spectral elements. In this work, we only focus on polynomial
non-conforming cases, but extension to geometric non-conforming cases is also
possible. Numerical validation is given in section~\ref{sec:tests}. Different
bases, quadratures, and shape types are tested to demonstrate that our
implementation is unified for general spectral elements. Finally, a detailed
performance benchmark of our initial matrix-free implementation is given to mark
our current progress.

\section{Basics and formulations}
\label{sec:basics}
\subsection{Spectral elements in Nektar++} \label{sec:Nektar}

The computing domain $\Omega_h$ is discretized into non-overlapping elements
$\Omega_e$: $\Omega_h = \bigcup \Omega_e$ with boundary $\partial
\Omega_h$. Similarly, element boundaries are denoted as $\partial
\Omega_e$. This section gives an overview of how Nektar++ unifies different
spectral elements. Details can be found in~\cite{karniadakis_spectralhp_1999}.

\subsubsection*{Bases and coefficient spaces}

The basis (or shape) functions $\phi(\boldsymbol{\xi})$ are defined on standard
elements $\Omega_{\text{st}}$ and the solution is approximated by expansions
\begin{equation} \label{eq:expansion}
  u_h(\boldsymbol{\xi},t)=\sum_{i}^{}{\phi_i(\boldsymbol{\xi})\hat{u}_i},
\end{equation}
where $\hat{u}_i$ is the coefficient related with the $i$-th basis. solutions are 
stored as a vector of coefficients for all bases and elements, ${\hat{\mathbf{u}}}$, 
known as \textit{coefficient} space. 

The construction of a general tensorial basis for a general element of dimension
$d$ is documented in~\cite{karniadakis_spectralhp_1999,laughton-2022}, and
involves a Duffy transformation~\cite{duffy-1982} between the reference element
$\Omega_{\text{st}}$ with coordinates $\boldsymbol{\xi}$, and the so-called {\em
  collapsed coordinate} system $\boldsymbol{\eta} \in [-1,1]^d$. Generally, this takes the form
\begin{equation}
  \phi_i(\boldsymbol{\xi}) = \psi_{i_1}(\eta_{1}) \psi_{i_2}(\eta_{2}) \cdots \psi_{i_d}(\eta_{d}),
\end{equation}
where $\psi_{i_k}(\eta_{k})$ is the basis function in coordinate direction $k$,
and for hexahedra and quadrilaterals, $\eta_{k}$ can be replaced by $\xi_k$. The
basis functions available in Nektar++ are related to element shape types and
summarised in Table~\ref{tab:elements & bases}. The modified $C^0$ \textit{modal
  basis} in~\cite{karniadakis_spectralhp_1999} were designed for CG methods, as
the basis functions possess a \textit{boundary-interior decomposition}, making
the system easier to assemble. In DG, this is not required but can be exploited
to simplify the trace contributions to the element.

\begin{table}
    \caption{Summary of elements and basis functions available in Nektar++ \cite{warburton_galerkin_1999}}
    \centering
    \begin{tabular}{llll}
    \toprule
        {\bf Shape type} & {\bf Basis type} & {\bf Tensor structure} & {\bf B/I decomposition} \\
   \midrule
        Quadrilateral/ & Modal        & standard & Yes \\
        Hexahedron     & Orthogonal   & standard & No  \\
                       & Lagrange-GLL & standard & Yes \\
                       & Lagrange-GL  & standard & No  \\
    \midrule
        Triangle/     & Modified modal  & generalized  & Yes \\
        Tetrahedron   & Orthogonal      & generalized   & No  \\
                      & Lagrange-Elec   & No            & Yes \\
    \midrule
        Prism         & Modified modal  & generalized & Yes \\
                      & Orthogonal      & generalized   & No  \\
                      & Lagrange-Elec   & No            & Yes \\
    \midrule
        Pyramid       & Modified modal  & generalized  & Yes \\
                      & Orthogonal      & generalized   & No  \\
    \bottomrule
    \end{tabular}
    \label{tab:elements & bases}
\end{table}

 \subsubsection*{Quadrature and physical spaces} 
 
\begin{figure}
    \centering
        \begin{subfigure}[b]{0.2\textwidth}
        \centering
        \includegraphics[width=\textwidth]{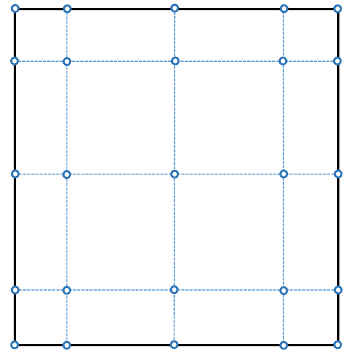}
        \caption{}
        \label{fig:GLL}
    \end{subfigure}
    \begin{subfigure}[b]{0.2\textwidth}
        \centering
        \includegraphics[width=\textwidth]{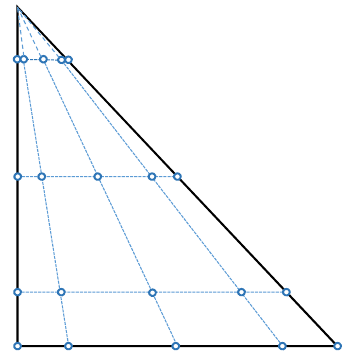}
        \caption{}
        \label{fig:GR}
    \end{subfigure}
    \begin{subfigure}[b]{0.2\textwidth}
        \centering
        \includegraphics[width=\textwidth]{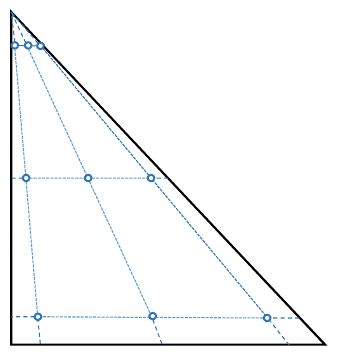}
        \caption{}
        \label{fig:GL}        
    \end{subfigure}
    \begin{subfigure}[b]{0.2\textwidth}
        \centering
        \includegraphics[width=\textwidth]{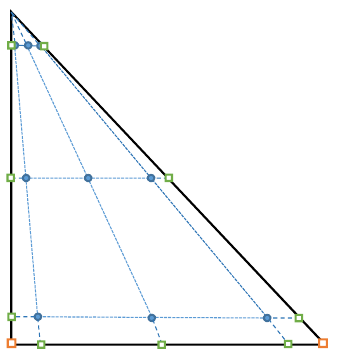}
        \caption{}
        \label{fig:GLwithEnd}
    \end{subfigure}
    \caption{The distribution of Gauss quadrature points for quadrilateral and triangle elements. (a) Gauss-Lobatto-Legendre (GLL) points in a quadrilateral; (b) Gauss-Radau (GR) points on dimension 2 in a triangle; (c) Gauss-Legendre (GL) points in a triangle; (d) Gauss points with additional endpoints for quadrature on both elements and boundaries, which still maintains a standard tensor structure.}
    \label{fig:GLGLLGR}
\end{figure}

In general, arithmetic operations must be done in \textit{physical} space within
${\mathbf{u}}$, a vector storing solution values at points
$\boldsymbol{\xi}_q$. These points are also responsible for integration within
the element, and one can increase the quadrature order to decrease aliasing
error. In Nektar++, multidimensional quadrature is constructed by a tensor
product of 1D quadrature
\begin{equation} \label{eq:iproduct}
  \int_{\Omega_e} u \, \text{d} \Omega \cong \sum_{p=1}^{N_{Q_1}} \sum_{q=1}^{N_{Q_2}} u(\xi_{1p}, \xi_{2q})\omega_{1p} \omega_{2q} J(\xi_{1p}, \xi_{2q})
\end{equation}
as shown in Fig.~\ref{fig:GLGLLGR}, where $N_{Q_i}$ are the number of quadrature
points in dimension $i$, $\omega_{ij}$ are the corresponding quadrature weights,
and $J$ is the determinant of Jacobian matrix
$|\frac{\partial (x_1,x_2)}{\partial (\xi_1,\xi_2)}|$. Gauss-Lobatto is
typically used, which includes endpoints and allows boundary conditions to be
easily imposed. However in non-tensorial elements, the Duffy transformation
leads to a removable singularity at collapsed vertices. We can avoid special
treatment of this point through the use of Gauss-Radau points, which exclude one
endpoint, or Gauss-Legendre points, which exclude both endpoints.

\subsection{DG discretisation of the Helmholtz equation using SIPG} \label{sec:SIPG}

\subsubsection*{Model problem and discretisation}
Consider the Helmholtz problem given by $ -\nabla^2 u + \lambda u= - f$ defined
on $\Omega$ with $\lambda \geq 0$. The discretisation by the symmetric interior
penalty Galerkin method (SIPG) on a single element reads~\cite{fehn_stability_2017}
\begin{multline}  \label{eq:IP}
    (\nabla{v_h},\nabla{u_h})_{\Omega_{e}}  + \lambda({v_h},{u_h})_{\Omega_{e}} \\
    -(\nabla{v_h},\frac{1}{2}[\![u_h]\!])_{\partial\Omega_e}
    -(v_h,\left\{\!\!\{ \nabla{u_h}\right\}\!\!\}\cdot{\boldsymbol{n}})_{\partial\Omega_e}
    +(v_h,\tau[\![u_h]\!]\cdot{\boldsymbol{n}})_{\partial\Omega_e}
  = -(v_h,f)_{\Omega_{e}}  , 
\end{multline}
where $v_h$ is the test function for $u_h$ space. The inner product is abbreviated by $(\boldsymbol{v}, \boldsymbol{u})_{\Omega}=\int_{\Omega}{ \boldsymbol{v} \cdot \boldsymbol{u} \text{ d}\Omega}$. 
In the context of DG, information on the current and its adjacent element is usually denoted by superscripts $^+$ and $^-$. 
The average operator $\{\!\!\{\cdot\}\!\!\}$ and jump operator $[\![\cdot]\!]$ are defined as 
\begin{equation}
    \{\!\!\{\boldsymbol{u}\}\!\!\}=\frac{\boldsymbol{u}^-+\boldsymbol{u}^+}{2}, \quad [\![\boldsymbol{u}]\!]=\boldsymbol{u}^-\otimes\boldsymbol{n}^-+\boldsymbol{u}^+\otimes\boldsymbol{n}^+  , 
\end{equation}
where $\boldsymbol{n}$ is the outward unit normal vector, and $\otimes$ denotes
the tensor product. On the boundaries, $\boldsymbol{u}^+$ and
$\nabla \boldsymbol{u}^+$ are determined by boundary conditions~\cite{fehn_stability_2017}.

Eq.~\eqref{eq:IP} is given from an element perspective. Summing over all
elements restores the primal formulation~\cite{arnold_unified_2002}

\begin{multline}  \label{eq:primal}
    (\nabla{v_h},\nabla{u_h})_{\Omega_h} + \lambda({v_h},{u_h})_{\Omega_h} \\
    - (\left\{\!\!\{ \nabla v_h \right\}\!\!\}, [\![ u_h ]\!])_{\Gamma} - ( [\![ v_h ]\!], \left\{\!\!\{ \nabla{u_h} \right\}\!\!\})_{\Gamma}  + \tau([\![ u_h ]\!], [\![ v_h ]\!])_{\Gamma} = -(v_h,f)_{\Omega_h}  , 
\end{multline}

where $\Gamma$ is the collection of all interior element interfaces and exterior
boundaries: the \textit{trace} space of the domain. To restore this bilinear
form, two overlapping interior element boundaries will be considered as the same
one; in other words, the inner product operator should be consistent for the
elements on the two sides to achieve a symmetric system. One primary reason for
choosing SIPG is that it can be solved by the conjugate gradient method, which
is less complicated to compute than other Krylov methods such as GMRES. However,
it is not guaranteed to converge if the system is not symmetric. 
This is in contrast to some other DG
applications, where an inaccurate flux evaluation may not ruin the solution
immediately. Our experience is that to achieve a truly symmetric SIPG
formulation is not trivial for non-conforming discretisations: a detailed
hands-on analysis is given in the appendix, showing how the flux evaluation
affects symmetry. There are two basic strategies to resolve this issue: either
to choose sufficient points to evaluate the quadrature accurately, or to unify
$\Gamma_l$ and $\Gamma_r$ and let the inner product be evaluated by the same
quadrature formula. The following section will discuss how we develop unified
interface flux evaluations based on these two strategies.

\section{Unified interface flux evaluation for general spectral elements} \label{sec:unified}

\subsection{From basic operators to optimised solver workflow}  \label{sec:workflow}

It is necessary to first present an overview of how a matrix-free DG solver is
implemented in Nektar++ before we discuss interface flux evaluation. 
The primary target is to evaluate the left-hand side (LHS) or right-hand
side (RHS) of a discretised equation. We prefer not to assemble the system
matrix in advance but evaluate the results via a series of basic operators
following its mathematical formula, which reduces memory access at runtime at
the cost of increasing floating point operations, thus increasing the arithmetic
intensity and the performance. To achieve this goal, the workflow must be carefully designed.
Nektar++ provides a series of core elemental operators:

\begin{itemize}
  \item \textbf{BwdTrans}: Evaluate the physical values $\boldsymbol{u}$ from
  the coefficients $\hat{\boldsymbol{u}}$ as in~\eqref{eq:expansion}.
  \item \textbf{PhysDeriv}: Compute the derivative of $\boldsymbol{u}$ using
  collocation differentiation~\cite{karniadakis_spectralhp_1999}. Note that the
  physical points always have a tensorial structure, so PhysDeriv is the same
  for any element shape.
  \item \textbf{IProduct}: Perform an inner product within this element
  according to~\eqref{eq:iproduct}. If we extract the quadrature metric $\omega$
  and $J$ out of the inner product, then this operation is essentially a
  transposed operation of BwdTrans.
\end{itemize}

\begin{figure}
    \centering
    \begin{subfigure}[b]{0.4\textwidth}
        \centering
        \includegraphics[width=\textwidth]{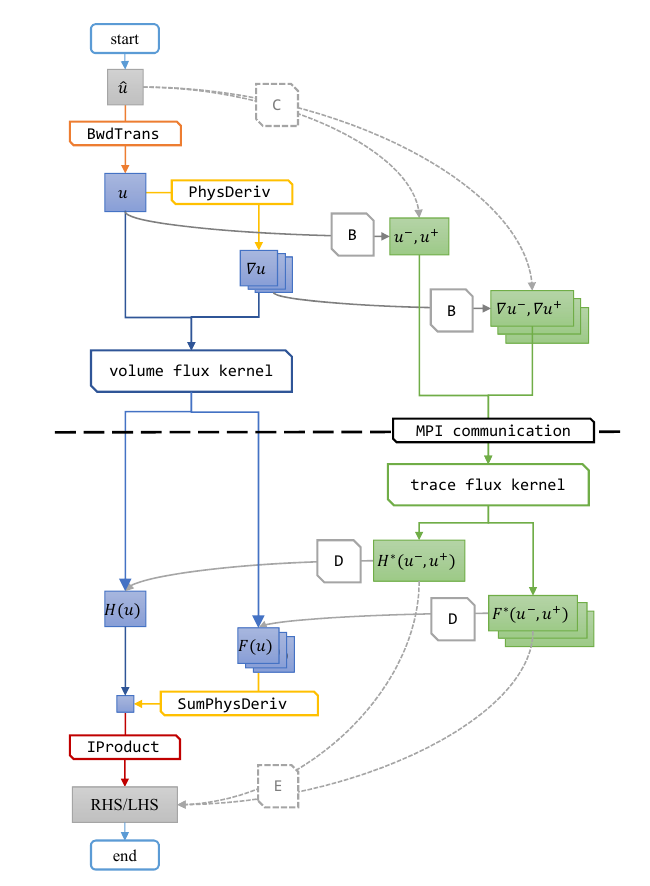}
        \caption{}
        \label{fig: essential workflow}
    \end{subfigure}
    \begin{subfigure}[b]{0.4\textwidth}
        \centering
        \includegraphics[width=\textwidth]{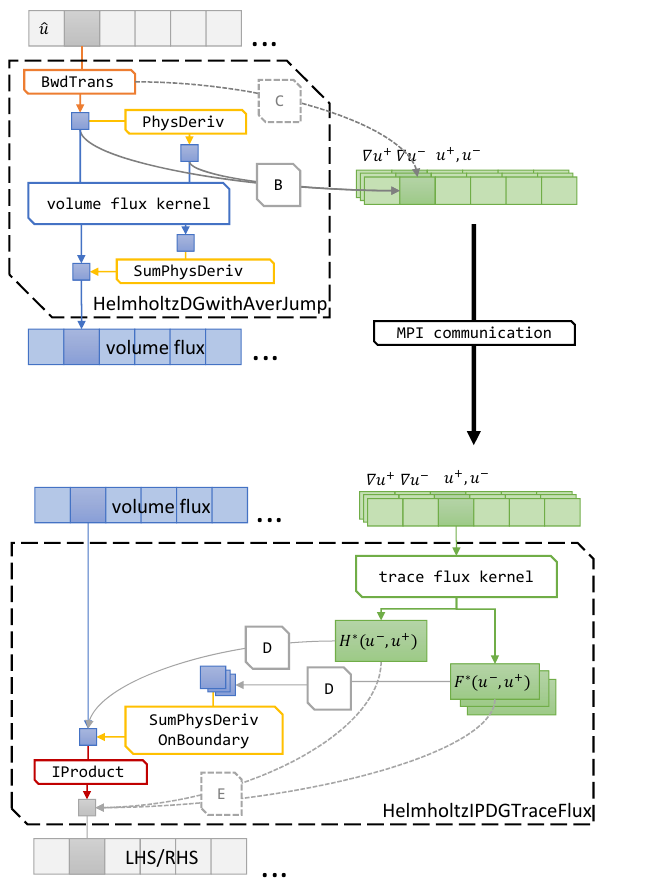}
        \caption{}
        \label{fig: actual workflow}
    \end{subfigure}
    \caption{The optimised workflow to evaluate the LHS or RHS in a DG formulation. (a) The basic design; (b) The actual design with fused kernel and cache blocking}
    \label{fig:workflow}
\end{figure}

Fig.~\ref{fig: essential workflow} shows the basic workflow to evaluate the LHS
or RHS of DG formulations; the actual implementation is shown in Fig.~\ref{fig:
  actual workflow} and will be discussed in the following section. Starting from
the solution unknowns $\hat{u}$, we call \texttt{BwdTrans} to get physical $u$
then \texttt{PhysDeriv} to get derivatives. Physical $u$ and $\nabla u$ are used
to evaluate the volumetric flux. Finally, we perform the inner product and add
the contribution to LHS or RHS, in the coefficient space. This workflow is not
unique: for example, one may consider combining \texttt{BwdTrans} and
\texttt{PhysDeriv} to get derivatives directly from coefficient space. However,
theoretical estimation~\cite{fischer_scalability_2020} and our practice show
this is more expensive than collocation differentiation \texttt{PhysDeriv} in
cases without over-integration.

As for interface flux term evaluation, we need additional operators \texttt{GathrInterp} and \texttt{ScatrInterp} (or \texttt{TracePhysEval} and \texttt{TraceIProduct}). Broadly speaking, a general interface flux evaluation takes three steps: 

\begin{enumerate}
    \item Gather trace physical data from element physical data (\texttt{GathrInterp}), or evaluate them directly from element coefficients (\texttt{TracePhysEval});
    \item Compute the flux on the trace physical space, varying according to the problems;
    \item Add trace physical data back to element trace physical space (\texttt{ScatrInterp}) and then perform inner product to get contributions to the LHS/RHS, or perform inner product on trace and transform the result back to element coefficient space directly (\texttt{TraceIProduct}).
\end{enumerate}
The relative merits of these operators will be discussed in the following sections.

\subsection{Obtaining trace data from elements}

\subsubsection*{Direct evaluation of trace physical solution \& derivatives}

\texttt{TracePhysEval} in Fig. \ref{fig: essential workflow} is essentially a special version of \texttt{BwdTrans} and \texttt{PhysDeriv}. If $\boldsymbol{\xi}_p$ covers all quadrature points on the trace and $\phi_i$ covers all bases of this element, then the resulting $\phi_i(\boldsymbol{\xi}_p)$ can be used to evaluate the physical solution on the trace. Sum-factorisation can also be applied to this operator:
\begin{equation}
    u(1,\xi_{2q},\xi_{3r}) = \sum_{i=1}^{N_{P}} \sum_{j=1}^{N_{P}} \sum_{k=1}^{N_{P}} \psi_i(1) \psi_j(\xi_{2q}) \psi_k(\xi_{3r}) \hat{u}_{ijk} \text{ , } 1\leq q,r \leq N_{Q_\Gamma} ,
\end{equation}
which evaluates all solution values on the face $\xi_1=1$, with a nominal
computing cost of $N_P^3+N_P^2N_{Q_\Gamma}+N_PN_{Q_\Gamma}^2$. 
\texttt{TraceIProduct} can be derived in the same way from~\eqref{eq:iproduct}.
As for the derivative on the trace, we need $\nabla \phi_j(\boldsymbol{\xi}_p)$ that covers all the trace quadrature points and all the element bases, which is equivalent to $\phi_i(\boldsymbol{\xi}_q)$ times $\nabla \ell_q(\boldsymbol{\xi}_p)$. 

\subsubsection*{Gather \& interpolate trace data from the element}

\texttt{GathrInterp} is available when the quadrature points include the endpoints, such as the Gauss-Lobatto points. In this case the boundary quadrature points are a subset of element quadrature points, so we can directly gather trace data from the elemental physical space without any additional computing cost. If the trace points are different from local surface points, then additional Lagrange interpolation is required. For example, to interpolate from $N_Q\times N_Q$ local face to $N_{Q_\Gamma} \times N_{Q_\Gamma}$ trace space, the nominal computing cost of interpolation by sum-factorisation is $N^2_Q N_{Q_\Gamma} + N_Q N_{Q_\Gamma}^2$. Recall that the trace space always has the same or higher quadrature order than its adjacent local element. The Lagrange interpolation will not deteriorate the polynomial order of the input, so we should obtain equivalent results as \texttt{TracePhysEval}. 

To use Gauss points in \texttt{GathrInterp} as well, one possible way is to add auxiliary endpoints to cover boundary quadrature, as shown in Fig. \ref{fig:GLwithEnd}. The resulting physical space still has a standard tensorial structure, so the existing sum-factorisation operators can be directly applied.


\subsection{Handling non-conforming interfaces}

\subsubsection*{Shared trace space} \label{sec:shared trace}

\begin{figure}
    \centering
    \includegraphics[width=0.7\textwidth]{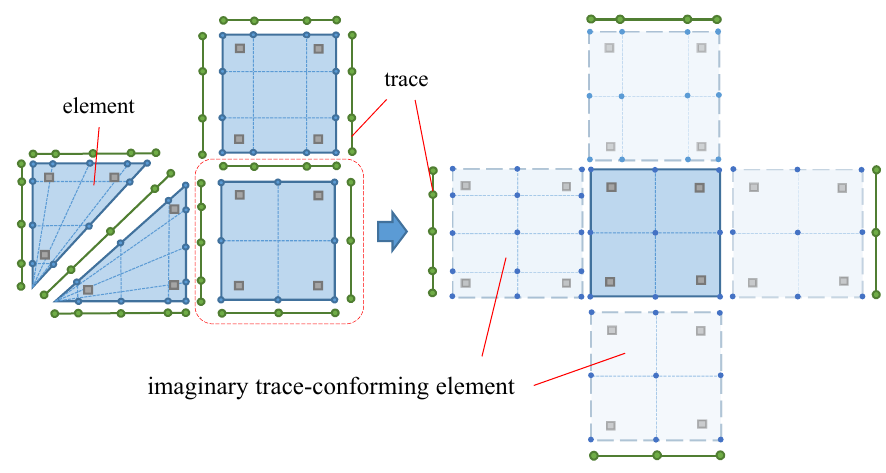}
    \caption{The concept of the shared trace space and trace-conforming element, which has the same bases and geometry as the local element, but the quadrature points are the same as the linked trace.}
    \label{fig:trace conform element}
\end{figure}

As discussed in section~\ref{sec:SIPG}, one strategy for a consistent flux evaluation on non-conforming interfaces is to unify the trace quadrature for two adjacent elements. So we begin by introducing the \textit{trace space} to Nektar++, a non-overlapping interface shared by two adjacent elements, as shown in Fig. \ref{fig:trace conform element}. In Nektar++, this concept is achieved by creating a list of trace \texttt{Expansion} objects, on which we can build required operators. As a general rule, the quadrature order on the trace should be at least the same as the higher-order side to achieve better accuracy.

\texttt{TraceIProduct} performs the inner product on the shared trace and projects results back to the local element coefficient space. Since the shared trace space is inconsistent with the local element on one side, we need to create a special set of \texttt{Expansion}, which may differ from the local elements, to get the corresponding $\phi_i(\boldsymbol{\xi}_q)$, and other geometric information used in operator \texttt{TraceIProduct} and \texttt{TracePhysEval}. They are called \textit{trace-conforming} element because the geometry and bases are the same as the local element, while the physical points are consistent with the trace. These \texttt{Expansion}s serve as a bridge between the shared trace and the local element.

\begin{figure}
    \centering
    \begin{subfigure}[b]{0.6\textwidth}
        \centering
        \includegraphics[width=\textwidth]{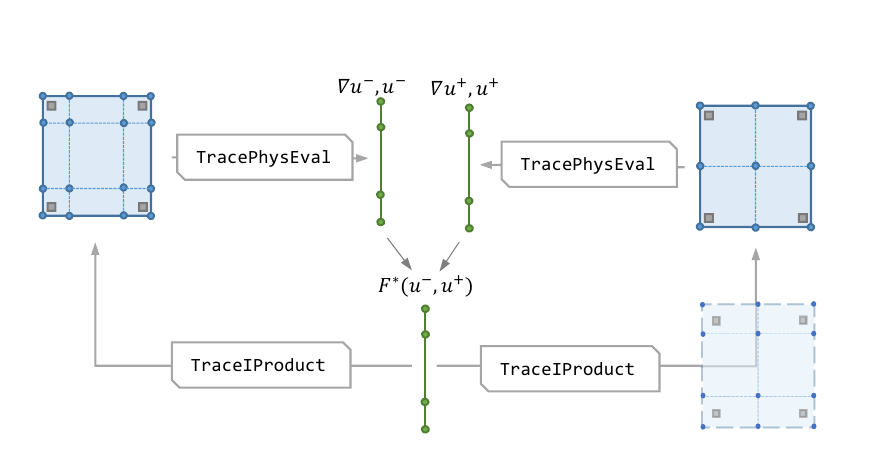}
        \caption{}
        \label{fig:shared trace}
    \end{subfigure}
    \begin{subfigure}[b]{0.6\textwidth}
        \centering
        \includegraphics[width=\textwidth]{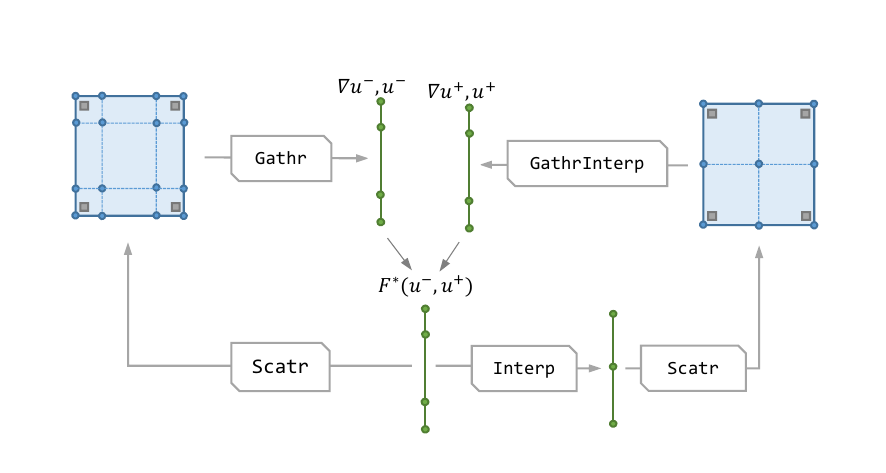}
        \caption{}
        \label{fig:p2p interp}
    \end{subfigure}
    \caption{Two different approaches handle non-conforming interfaces. (a) Flux evaluation via shared trace; (b) Flux evaluation by point-to-point interpolation.}
    \label{fig:two approaches to handle non-conforming interfaces}
\end{figure}

The shared trace space is conceptually the same as the \textit{mortar elements},
which are traditionally used to bridge between non-conforming interfaces in
continuous spectral element meshes. In the mortar element method, we first
project solutions from the local elements on the two sides to the mortar
element. Then we evaluate the flux in the mortar space and finally project the
flux back to the element of each side. For CG or DG methods, an $L^2$ projection
is a natural choice~\cite{kopriva_computation_2002} and is given by

\begin{equation} \label{eq:mortar project}
    \hat{\boldsymbol{u}}^{\Xi} = (\mathbf{M}^{\Xi})^{-1} \mathbf{S}^{\Omega_e \rightarrow \Xi} \hat{\boldsymbol{u}}
\end{equation}

where the superscript $\Xi$ denotes the objects of the mortar element,
$\mathbf{M}^{\Xi}$ is the standard elemental mass matrix of mortar element
${M}^{\Xi}_{ij} = (\phi^{\Xi}_i,\phi^{\Xi}_j)_{\Xi}$, and
$\mathbf{S}^{\Omega_e \rightarrow \Xi}$ is the transfer mass matrix from local
element to mortar element $S^{\Omega_e \rightarrow \Xi}_{ij} = (\phi^{\Xi}_i,\phi_j)_{\Xi}$. 
Previous studies considered nodal elements with collocation, so that the mass matrix is
diagonal and easy to apply. However, in a general spectral element framework,
the mass matrix is dense and generally differs per element, making computing and
multiplying the inverse mass matrix in every projection an expensive part of the
method.
Additionally, as mentioned in section~\ref{sec:Nektar}, the flux calculation
needs to be performed in physical space rather than coefficient space, which
leads to an additional call to \texttt{BwdTrans}. Typically, we only need to
perform the inner product at the final step and multiply by the inverse mass
matrix if required. The projection between the coefficient space is sub-optimal
and should be avoided if possible. To achieve this, Eq.(~\ref{eq:mortar project})
can be rewritten as
\begin{equation}
    \int_{\Xi} \phi^{\Xi}_i(\boldsymbol{z}) u^{\Xi}(\boldsymbol{z}) {\rm d}\Xi = \int_{\Xi} \phi^{\Xi}_i(\boldsymbol{z}) u(\boldsymbol{\xi}) {\rm d}\Xi
\end{equation}

In the case that the mortar element is aligned with the local element, the local 
coordinate $\boldsymbol{z}$ and $\boldsymbol{\xi}$ are the same. 
If our goal is to retrieve the physical solution on the mortar elements 
that satisfy the above $L^2$ relation, we can imprint the mortar quadrature points 
$\boldsymbol{z}_q$ to the local element and evaluate the solution on those points:

\begin{equation}
    u^{\Xi}(\boldsymbol{z}_q) = u(\boldsymbol{z}_q) = \sum_{j} \phi_j(\boldsymbol{z}_q) \hat{u}_j.
\end{equation}

This is the approach we adopt in our current design, so that the flux evaluation
via shared trace space in the present work is mathematically equivalent to the
mortar element method, but simplified in computation. It can also be extended to
handle geometric non-conforming interfaces as the mortar element does, although
this is beyond the scope of the current work.

\subsection*{Point-to-point interpolation} \label{sec:p2pinterp}

In \texttt{ScatrInterp}, we first use Lagrange interpolation to transform data from
the adjacent trace space into the local trace space, and then scatter the trace data
back to the local element, with inner products evaluated in the local element
space separately, as shown in Fig. \ref{fig:p2p interp}. The drawback is that we
cannot achieve symmetrical systems in all cases, due to the interpolation from higher order to lower order. 
Following the analysis in the appendix, the matrix $\mathbf{M}_C$ and $\mathbf{M}'_C$ in the case of
point-to-point interpolation becomes

\begin{align*}
  M^T_{C(j,i)}=\sum^{N_{Q_l}}_{p=0} \phi^{l}_{j} (\xi_{p})  (\omega_{p} J_{p}) \boldsymbol{n}^{l}_{p} \cdot \sum^{N_{Q_r}}_{q=0} \ell_{q}(\xi_{p}) \nabla \phi^{r}_{i} (\xi_{q}) \\
  M'_{C(i,j)}=\sum^{N_{Q_r}}_{q=0} \nabla \phi^{r}_{i} (\xi_{q})  \cdot  (-\boldsymbol{n}^{l}_{q}) (\omega_{q} J_{q}) \sum^{N_{Q_l}}_{p=0} \ell_{p}(\xi_{q})  \phi^{l}_{j} (\xi_{p}) 
\end{align*}

To ensure $\mathbf{M}^T_{C} = -\mathbf{M}'_{C}$, we need to satisfy two conditions:
\begin{enumerate}
    \item $N_{P_l} < N_{P_r} \leq N_{Q_l} < N_{Q_r}$ or $N_{P_l} > N_{P_r} \geq N_{Q_l} > N_{Q_r}$. This is obvious because the polynomial degree should not decrease after interpolation from the higher-order to the lower-order side.
    \item The quadrature on the lower-order local trace should still be accurate for the integrand. The integrand can be up to $(N_{P_r}-1)+(N_{P_l}-1)+3p_{\text{geom}}$, where $p_{\text{geom}}$ is the polynomial degree of geometric mappings, and the multiple 3 includes $J$, $\boldsymbol{n}$ and $\frac{\partial \boldsymbol{x}}{\partial \boldsymbol{\xi}}$. The quadrature degree of exactness should be no less than that.
\end{enumerate}

\begin{figure}
    \centering
    \includegraphics[width=0.6\textwidth]{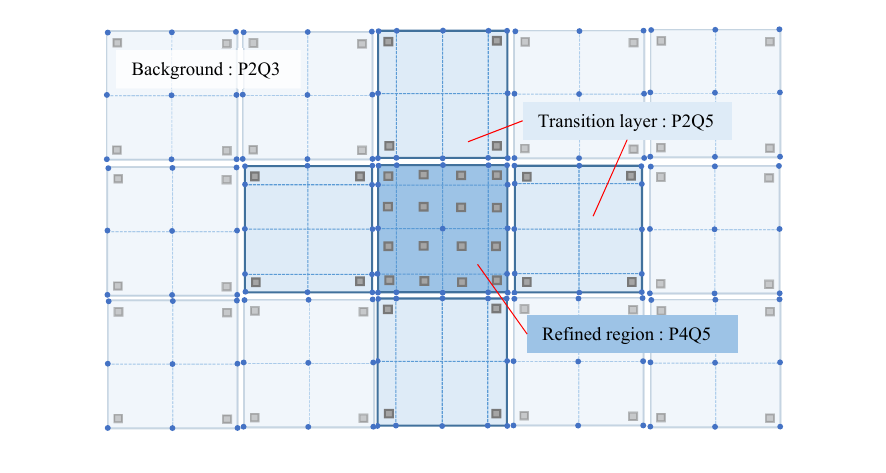}
    \caption{The transition layer around a $p$-refined region. The elements in the transition layer have the same points as the $p$-refined region, but the same coefficients as the background region.}
    \label{fig:transition layer}
\end{figure}

These conditions can be difficult to fulfill with arbitrary element shapes and
orders. However, the reason we still consider this method is that it is simple
to implement and can be easily applied to arbitrary non-conforming interfaces,
particularly in 3D where mortar elements of unstructured triangles are very hard to construct. 
Moreover, this approach is also easier to achieve higher performance, as shown in section \ref{sec:matfree}.
A practical tip we can suggest to help alleviate the symmetry error is to insert a 
\textit{transition layer}, as seen in Fig. \ref{fig:transition layer},  
at a price of greater computing cost, which is generally sufficient to resolve
conjugate gradient convergence issues in the current application. 
For example, a $P4Q5$ (4 coefficients, 5 GLL points) element 
cannot be the neighbour to the $P2Q3$ elements. After adding
such a transition layer, which has the same coefficients as the lower-order
elements but the same quadrature points as the higher-order element, the direct
neighbours now become $P2Q5$-$P4Q5$ and $P2Q5$-$P2Q3$. The former has a
conformal trace space. As for the latter, we only need to ensure the $P2Q3$ has
sufficient quadrature accuracy for its local trace evaluation to achieve a
symmetric system. This makes sense in adaptive $p$-refinement since the
background mesh $P2Q3$ will not be affected by the local refinement.

\subsection{Other implementation issues}

\subsubsection*{Singular vertices in collapsed coordinates} 
\begin{figure}
    \centering
    \begin{subfigure}[b]{0.25\textwidth}
        \centering
        \includegraphics[width=\textwidth]{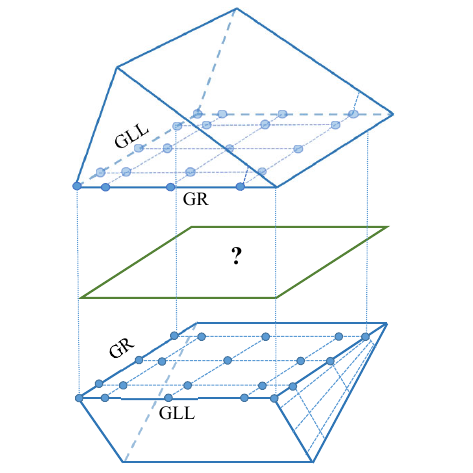}
        \caption{}
        \label{fig:GR in prism trace}
    \end{subfigure}
    \begin{subfigure}[b]{0.25\textwidth}
        \centering
        \includegraphics[width=\textwidth]{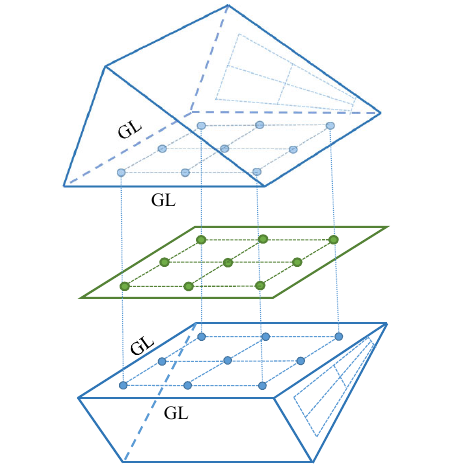}
        \caption{}
        \label{fig:GL in prism trace}
    \end{subfigure}
    \caption{If collapsed coordinates are not aligned to the adjacent element, we
        cannot construct trace points that conform to Gauss-Radau points on both
        sides as shown in (a), but this is possible with rotationally-symmetric Gauss
        points, as shown in (b).}
\end{figure}

As mentioned
previously, we avoid evaluating the derivative on singular vertices in collapsed
coordinates of a simplex using either Gauss or Gauss-Radau points. 
Generally speaking, when constructing the trace-conforming elements,
if the point distribution on a face or an edge is not rotational symmetric, 
then the local points on the two sides of the shared trace may be misaligned.
Nektar++ can adjust the vertex of tetrahedra to ensure the top singular vertex of
adjacent elements is always aligned. However, this cannot cover all the cases.
Fig. \ref{fig:GR in prism trace} reveals an exception that
that we cannot construct trace quadrature points that conform to both
sides if Gauss-Radau points are used on a quadrilateral face due to non-symmetrical point
distribution. The best and cheapest solution is using Gauss points, which are
rotationally symmetric, as shown in Fig. \ref{fig:GL in prism
  trace}.

\subsubsection*{Exploiting boundary/interior decomposition}
Boundary/interior decomposition is not required for DG. But with this property, the solution on the boundary only depends on certain coefficients that correspond to that boundary surface, known as the boundary modes \cite{warburton_galerkin_1999}. This immediately provides an opportunity to reduce \texttt{BwdTrans} computing cost in the \texttt{TracePhysEval}. We can use a smaller basis matrix $\phi_i(\boldsymbol{\xi}_q)$ where $\phi_i$ only includes boundary bases. The matrix size is reduced to $N^{d-1}_P \times N^{d-1}_{Q_\Gamma}$ and the sum-factorisation cost is $N_P N_{Q_\Gamma}^2 + N_P^2 N_{Q_\Gamma}$ in a hexahedron. Such boundary/interior decomposition is only valid for solution variables, but not for derivatives. So the overall cost reduction is modest.

\section{Numerical validation} \label{sec:tests}   

In this section, we provide a numerical validation of the concepts introduced in
the previous section by examining convergence order for the Helmholtz
equation. To achieve this, we consider the sinusoidal solution
\begin{equation} \label{eq:sinxsinysinz}
    u(x,y,z) = \sin(k x)\sin(k y)\sin(k z).
\end{equation}
with a corresponding manufactured forcing term $f$ given as
\begin{equation} \label{eq:force-sinxsinysinz}
    f = -(\lambda + 3k^2)\sin(k x)\sin(k y)\sin(k z).
\end{equation}

\subsection{Validation and comparison with point-to-point interpolation}

\begin{figure}
    \centering
    \begin{subfigure}[b]{0.32\textwidth}
        \centering
        \includegraphics[width=\textwidth]{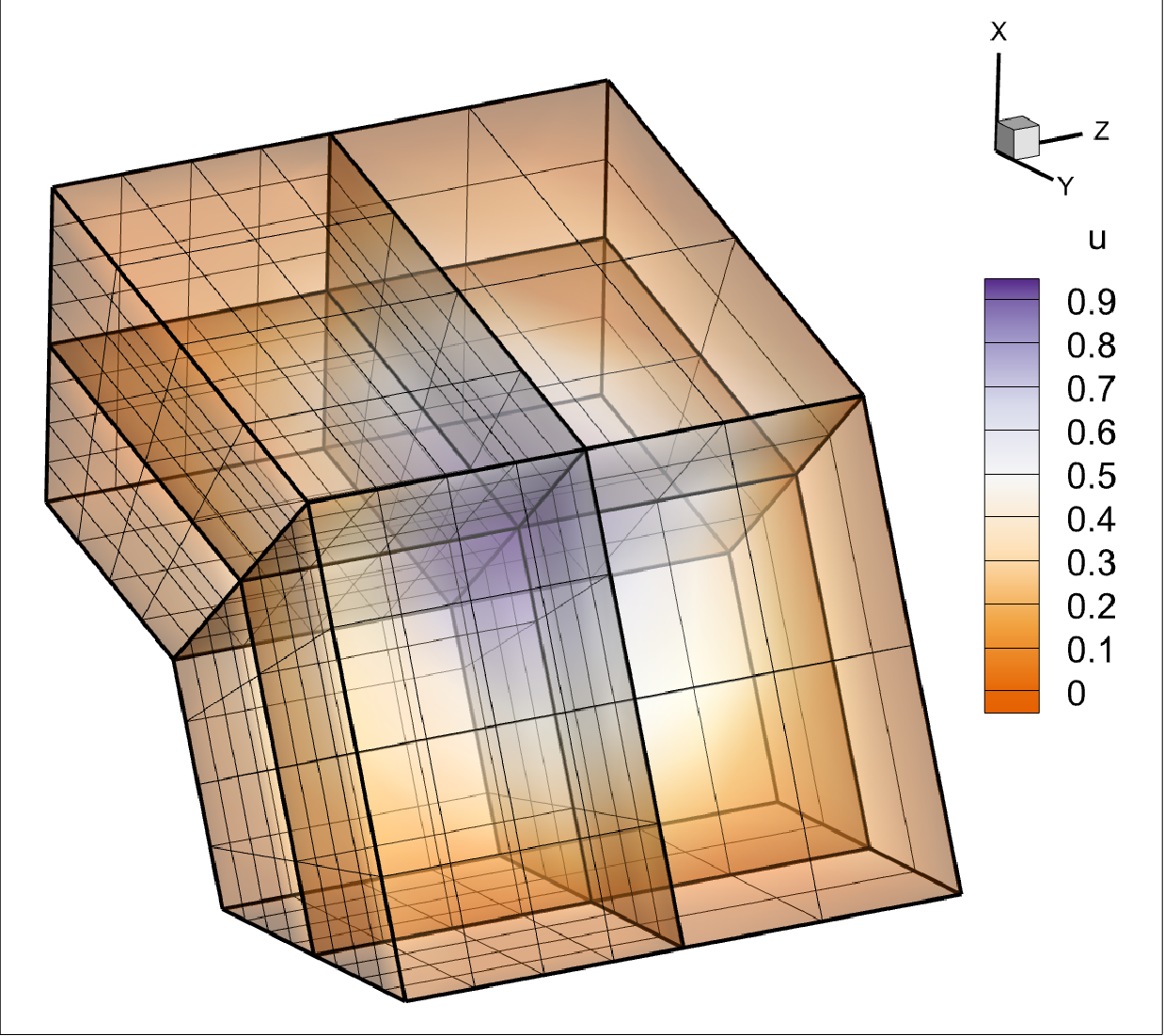}
        \caption{}
    \end{subfigure}
    \begin{subfigure}[b]{0.315\textwidth}
        \centering
        \includegraphics[width=\textwidth]{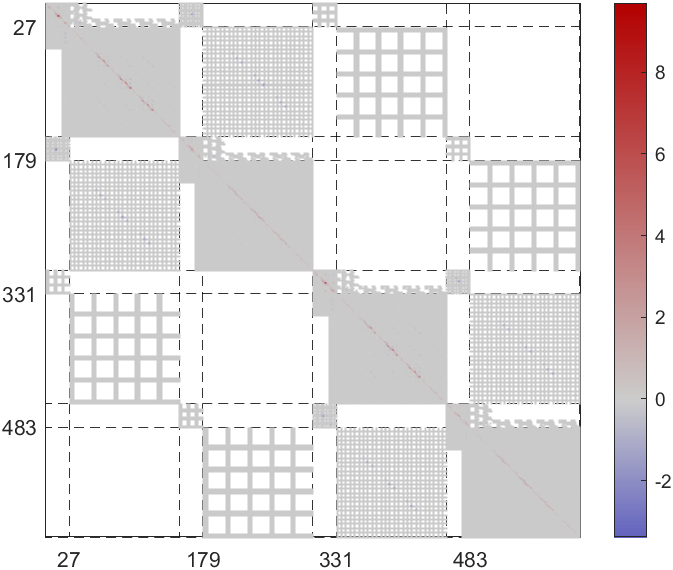}
        \caption{}
    \end{subfigure}
    \begin{subfigure}[b]{0.325\textwidth}
        \centering
        \includegraphics[width=\textwidth]{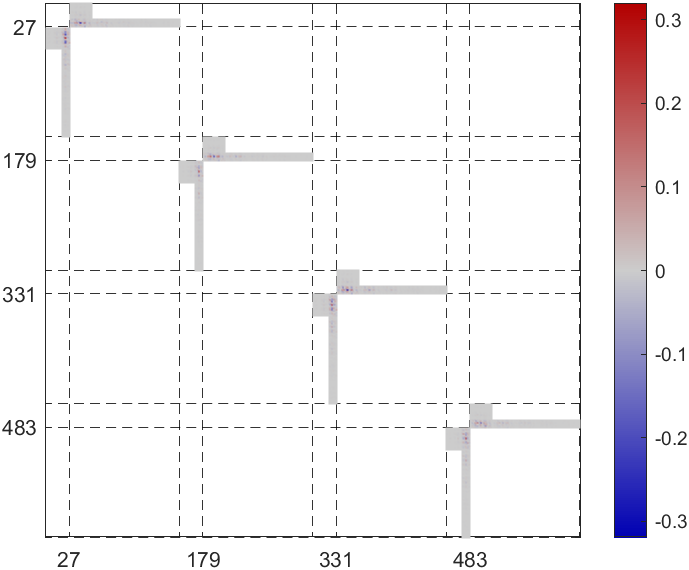}
        \caption{}
    \end{subfigure}
    \caption{A critical case for point-to-point interpolation: $P3Q4$ mixed with $P5Q6$. (a) The computing domain; (b) Sparsity pattern of the system matrix $\mathbf{A}$; (c) Sparsity pattern of $\mathbf{A}-\mathbf{A}^T$}
    \label{fig:system matrix}
\end{figure}

We start with basic validation by checking if the system matrix is symmetric in
the cases of local $p$-refinement with $k=\pi/2$. The domain is a sector filled
with $2^3$ linear-shape hexahedra, where the elements use a Lagrange basis and
independent GLL quadrature points. Boundary conditions are set to Dirichlet
conditions consistent with the manufactured solution. To test the non-conforming
implementation, half of the domain has a different polynomial order. For
point-to-point interpolation, the case $P3Q4$-$P5Q6$ is critical, which leads to
a non-symmetric system, as shown in Fig.~\ref{fig:system matrix}. The
non-symmetric entries all come from the off-diagonal part, such as the matrices
$M_C$ $M_D$ and $M_G$ in the appendix. Although the difference in
$\|\mathbf{A} - \mathbf{A}^T\|_1$ is less than $3\%$ of $\|\mathbf{A}\|_1$, it
is enough to prevent convergence, although solution via GMRES at equivalent
accuracy is possible, as shown in Table~\ref{tab:compare p2p interp}. The table
also shows that the system is exactly symmetric with $P3Q5$, which follows our
previous analysis.

\begin{table}
 \caption{The analysis results of mixed-order cases shown in Fig. \ref{fig:system matrix}, by point-to-point interpolation. All cases use GLL points and Lagrange bases.}
  \centering
  \begin{tabular}{llll}
    \toprule
    coefficient and quadrature & $P3Q4$-$P5Q6$   & $P3Q5$-$P5Q6$    & $P4Q5$-$P5Q6$ \\
    \midrule
    $N_{P_r} \leq N_{Q_l}$ ?    &   No            & Yes          & Yes        \\
    required quadrature degree  & $2+4+1=7$       & $2+4+1=7$    & $3+4+1=8$   \\
    actual quadrature degree   & $4\times2-3=5$  & $5\times2-3=7$ & $5\times2-3=7$ \\
    $\frac{\|\mathbf{A} - \mathbf{A}^T\|_1}{\|\mathbf{A}\|_1}$     &   $0.0252$    &   $2.5\times10^{-18}$   &  $0.0013$   \\
    CG iterations to $10^{-9}$   & -              & 62            & 69              \\
    GMRES iterations to $10^{-9}$ & 61             & 62            & 63              \\
    $L_2$ error                  & $1.52\times10^{-2}$ & $1.08\times10^{-2}$ & $8.96\times10^{-4}$ \\
    \bottomrule
  \end{tabular}
  \label{tab:compare p2p interp}
\end{table}

\subsection{$p$-refinement convergence rate}

We now consider a convergence study of the SIPG Helmholtz solver with local
polynomial refinement, for which we select a larger wavenumber $k=2\pi$, so as
to require higher resolution to attain accurate results. The cube domain is
filled with $N_x^3$ structured hexahedra. We subdivide each hexahedron into
either 6 tetrahedra or 2 prisms to create meshes with different element shapes
but similar resolutions. All elements use the orthogonal basis, GL points and
$N_Q=N_P$. Half of the domain is refined by adding 1 to the element order.

\begin{figure}
    \centering
    \begin{subfigure}[b]{0.32\textwidth}
        \centering
        \includegraphics[width=\textwidth]{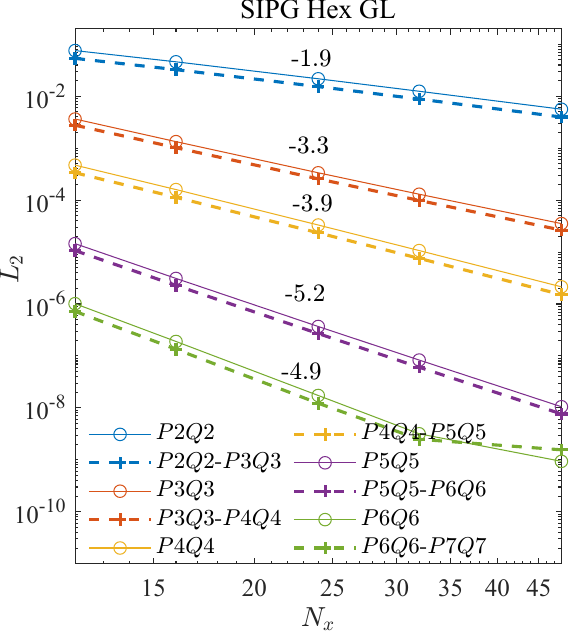}
        \caption{}
    \end{subfigure}
    \begin{subfigure}[b]{0.32\textwidth}
        \centering
        \includegraphics[width=\textwidth]{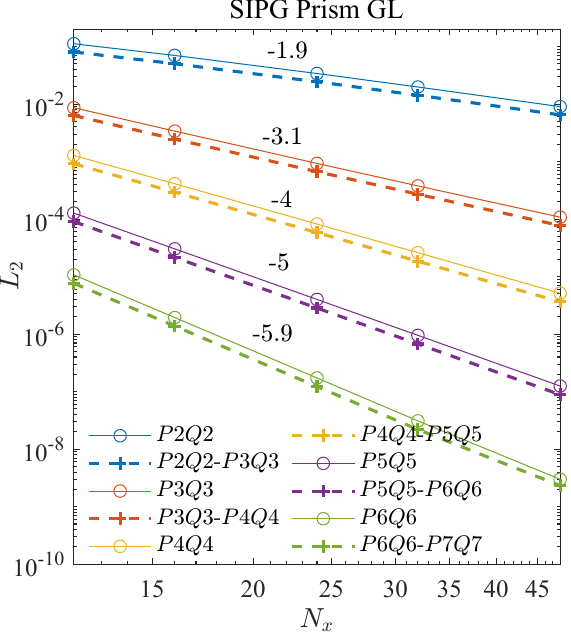}
        \caption{}
    \end{subfigure}
    \begin{subfigure}[b]{0.32\textwidth}
        \centering
        \includegraphics[width=\textwidth]{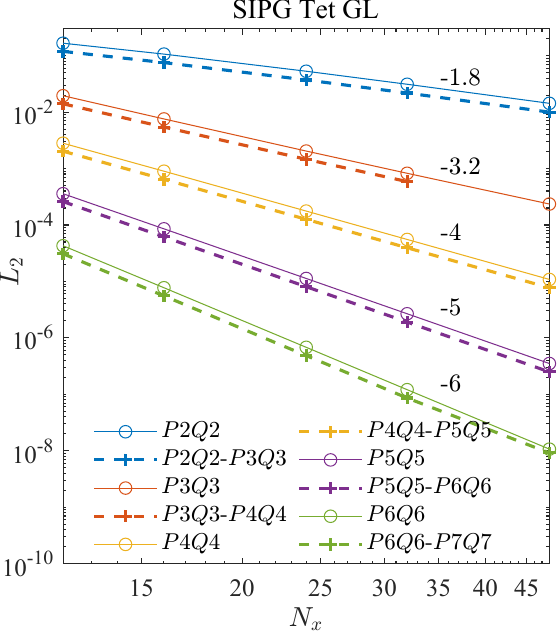}
        \caption{}
    \end{subfigure}
    \caption{The convergence rate of a sinusoidal problem with local $p$-refinement. (a) Hexahedral mesh; (b) Prismatic mesh  (c) tetrahedral mesh. The numbers are the slopes of the curves.}
    \label{fig:convergence-sinx}
\end{figure}

Fig.~\ref{fig:convergence-sinx} presents the relation between $L^2$ error and
mesh resolution $N_x$ with different polynomial orders and element shapes. From
the slopes given in the figure, the SIPG Helmholtz solver achieves around $N_P$
order convergence rate as expected. The locally refined cases also show a
similar convergence rate with respect to the non-refined cases, and the error is
only slightly lower, which is expected since only half the domain is
refined. Typically the error contributed by higher-order discretisation is
significantly smaller than the lower-order one, so the total error is dominated by
the lower-order region.

\begin{figure}
    \centering
    \begin{subfigure}[b]{0.32\textwidth}
        \centering
        \includegraphics[width=\textwidth]{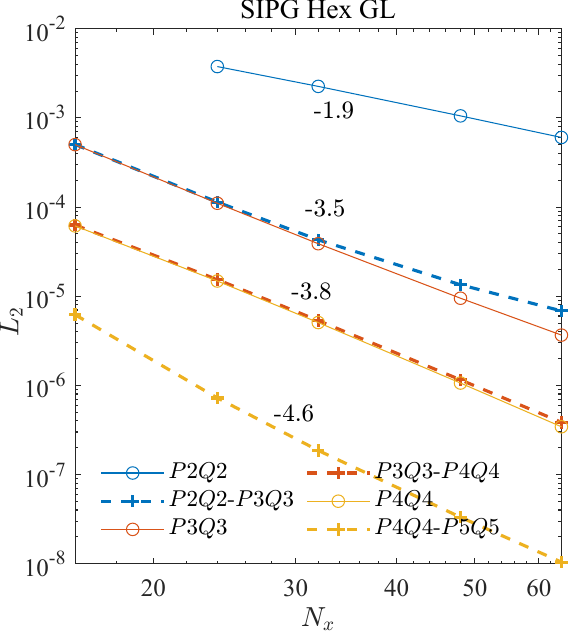}
        \caption{}
    \end{subfigure}
    \begin{subfigure}[b]{0.32\textwidth}
        \centering
        \includegraphics[width=\textwidth]{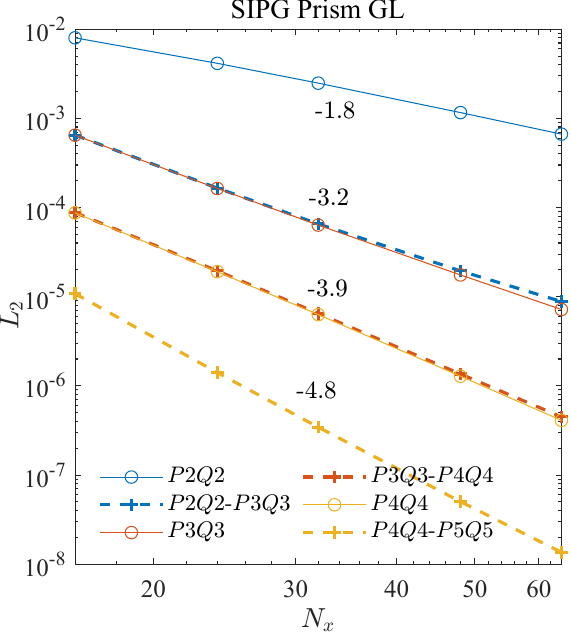}
        \caption{}
    \end{subfigure}
    \begin{subfigure}[b]{0.32\textwidth}
        \centering
        \includegraphics[width=\textwidth]{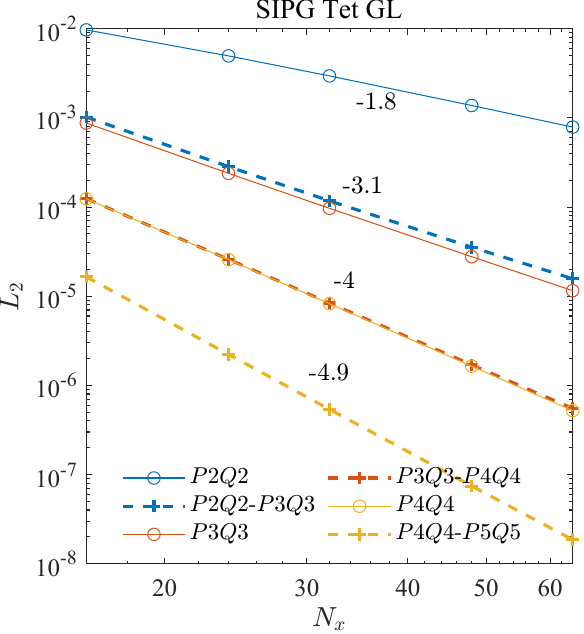}
        \caption{}
    \end{subfigure}
    \caption{The convergence rate of the Gaussian pulse problem with local $p$-refinement. (a) Hexahedral mesh; (b) Prismatic mesh;  (c) tetrahedral mesh. The numbers are the slopes of the curves.}
    \label{fig:convergence-Gaussian}
\end{figure}

To better show the effectiveness of local $p$-refinement, we manufacture a problem with a Gaussian pulse solution
\begin{equation}
    u(x,y,z)=\text{exp}\left( \frac{x^2+y^2+z^2}{a^2} \right)
\end{equation}
within the same cube domain $[-1,1]^3$. We select a character radius $a=0.2$ and the refined 
region is a cube in the middle of the domain $[-0.5,0.5]^3$, enough to contain the pulse. 
Outside the refined region, the solution is nearly zero, so the error mostly depends on the refined region.
This is confirmed in Fig.~\ref{fig:convergence-Gaussian}, where the locally refined case now achieves error magnitudes 
and convergence rates similar to the uniform high-order case, but with a lower number of degrees of freedom.

\section{Matrix-free implementation and performance benchmark} \label{sec:matfree}

In this section, we examine the high-performance implementation of the
formulation of the previous sections, and outline the results from benchmarking that
show the effective performance of the method on unstructured meshes.

\subsection{Design choices in Nektar++}

Efficient matrix-free kernels for core elemental operations are the foundation to achieve high performance, and Nektar++ already has these operators for various shape types \cite{moxey_efficient_2020}. The design choices we made in Nektar++ for any matrix-free implementations are:
\begin{itemize}
    \item Disassemble a complex operator to lower-level operators like \texttt{BwdTrans}, where sum-factorisation or collocation can be applied, and basis or geometric data can be shared across elements or operators.
    \item Evaluate all other complex coefficients based on the pre-computed shared geometric information, Jacobian $J$ and derived factor $\frac{\partial \boldsymbol{x}}{\partial \boldsymbol{\xi}}$. This is a trade-off between data size and computing cost. 
    \item Provide specialized kernels for \textit{regular} and \textit{deformed}
    elements, since geometric information is constant for regular shape but spatially varying in deformed (curvilinear) elements.
    In the former case, we only need to load one Jacobian value or one face normal, and apply it to all points, which
    significantly reduces the memory load demands and increases arithmetic intensity.
    \item Elements of the same shape, basis and order are grouped and processed
    together, so that we can reuse the basis data when processing each element. To
    reuse geometric data as well and avoid writing large intermediate results to
    main memory, we try to perform as many operations as possible on a small group
    of elements before we write final results to memory and move to the next. This
    is illustrated in Fig. \ref{fig: actual workflow} and called cache blocking or
    kernel fusing in some literature.  Kernel functions are force-inlined to
    reduce overheads.
    \item Interleave data layout across elements for more efficient SIMD vectorisation. 
    For example, the data with the same local index from 4 elements are stored contiguously in
    memory. These data will be directly loaded into vector registers and processed
    by SIMD (single instruction, multiple data) instructions such as AVX2 or AVX512 
    (Advanced Vector eXtensions). In other words, each operator in Nektar++ can process 4 to 8 elements at a time.
    \item Generate separate kernels for different orders so that the loop bounds and local temporary memory sizes are all compile-time constants. In these cases, compilers can more aggressively optimise through loop unrolling and fusion.
\end{itemize}

As for parallel execution, Nektar++ only exchanges trace data with adjacent mesh
partitions, instead of setting up \textit{ghost} elements around the partition
and exchanging the data of a whole element. This reduces the data size in
communications, especially for high orders. Synchronisation is required during
communication to avoid data racing. (We note it is also possible to overlap
communication and computation to hide the communication costs, although we do
not do this presently.) Our full workflow (\texttt{LhsEval}) is split into three
sequential parts as shown in Fig.~\ref{fig: actual workflow}:
\begin{enumerate}
  \item \texttt{HelmholtzDGwithAverJump}: Evaluates Helmholtz volume flux and gets the average or the jump data on the traces.
  \item \texttt{MPIexchange}: Sends trace data to adjacent partitions and also receives data from them.
  \item \texttt{HelmholtzIPDGTraceFlux}: Evaluates all the interface fluxes and add their contribution back to the element space.
\end{enumerate}
Both \texttt{HelmholtzDGwithAverJump} and \texttt{HelmholtzIPDGTraceFlux} are optimised fused kernels that process only a small
group of elements at a time. 
We minimise the access to the global memory to only inputs and outputs at the beginning or the end of the fused kernel. 
The sub-operators only write intermediate results to a small temporary memory, which should stay in the cache and be accessed efficiently by the current or the next sub-operator.
In this particular SIPG Helmholtz solver, we only need $\left\{\!\!\{ \nabla u \right\}\!\!\} \cdot \boldsymbol{n}$ and $[\![u]\!] \cdot \boldsymbol{n}$ to evaluate all the flux shown in Eq. \eqref{eq:IP}. So instead of exchanging four primitive variables (one $u$ and three $\nabla u$), we only exchange two variables to save communication costs.

We have two options to obtain trace data (\texttt{GathrInterp},
\texttt{TracePhysEval}) and two options to handle non-conforming interfaces
(\texttt{ScatrInterp}, \texttt{TraceIProduct}).  As for implementation, the
concerns focus on ease of integration into the framework and the potential to
achieve higher performance.  The first observation is that when interfaces are
conforming, \texttt{GathrInterp} and \texttt{ScatrInterp} only involve memory
transfer and have no floating-point cost. After \texttt{ScatrInterp}, only 1
\texttt{IProduct} on the whole element is performed. In contrast, we have to do
\texttt{TraceIProduct} for every edge or face of the element, which is more
expensive, as confirmed in both previous studies~\cite{laughton_comparison_2021}
and the present study. We therefore choose \texttt{GathrInterp} and
\texttt{ScatrInterp} as the primary design, with \texttt{TraceIProduct} only
enabled on non-conforming interfaces if required.

\begin{figure}
    \centering
    \begin{subfigure}[b]{0.41\textwidth}
        \includegraphics[width=\textwidth]{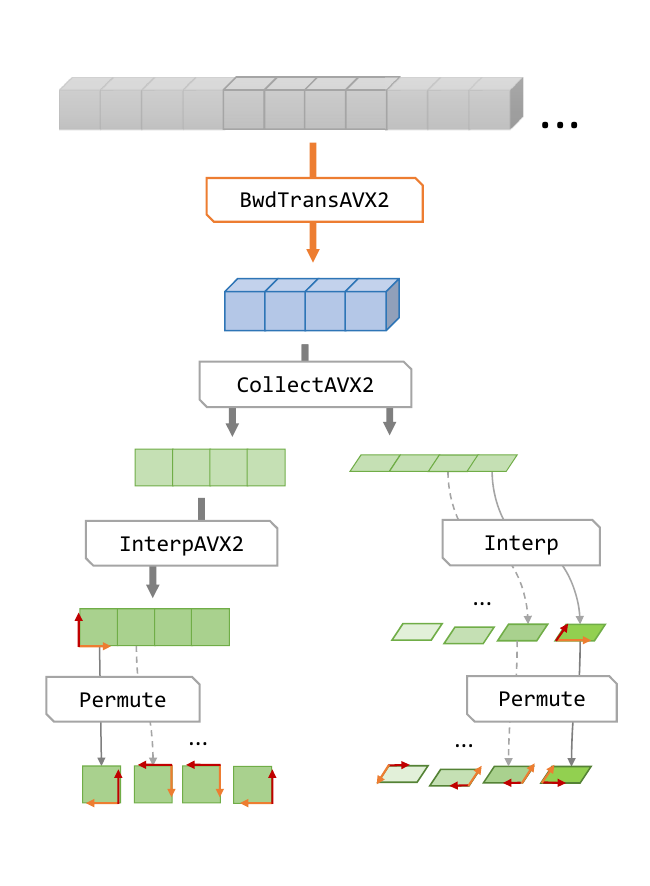}
        \caption{}
        \label{fig:vector strategy operB}
    \end{subfigure}
    \begin{subfigure}[b]{0.33\textwidth}
        \includegraphics[width=\textwidth]{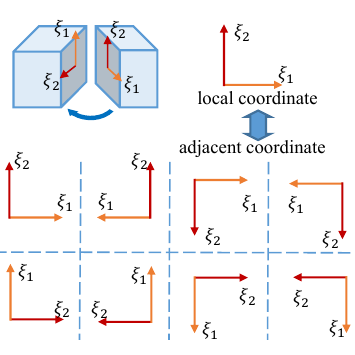}
        \caption{}
        \label{fig:orientations}
    \end{subfigure}
    \caption{(a) The vectorisation strategy for gathering and interpolating trace data, and (b) the orientations between adjacent faces.}
\end{figure}

Our matrix-free implementation of \texttt{GathrInterp} and \texttt{ScatrInterp}
follows the design choices above. In an arbitrary unstructured mesh, two
adjacent elements may have different trace spaces and also different
orientations, as shown in Fig.~\ref{fig:orientations}. In three-dimensional
space, there can be at most 8 different orientations, but usually, only a few
orientations may appear in a practical mesh, and can be further reduced by
deliberately swapping element local coordinates during the pre-processing
stage. So the detailed process of operator B actually takes three steps:

\begin{enumerate}
  \item {\bf Gather:} Gather local trace physical data from the element's
  physical space.
  \item {\bf Interpolation:} If local and adjacent trace spaces differ,
  perform interpolation.
  \item {\bf Permutation:} If the local trace has different orientations,
  permute the trace data so that the storage order is the same as the adjacent
  space.
\end{enumerate}

If the $i$-th trace of all elements in this group has the same interpolation configuration, we can perform gathering and interpolation in batch mode, using explicit SIMD instructions. Otherwise, only task 1 can be processed in batch mode and the other two must be done on a single trace at a time. This is illustrated in the Fig. \ref{fig:vector strategy operB}. The runtime efficiency, therefore, depends on how many elements of the same type and how many traces of the same type can be grouped together, which should be done during the preprocessing of the mesh. 

Note that we do not claim the above design choices are the best overall: indeed,
several were made to maintain compatibility with the existing codebase. The
efficient implementation of DG operators, especially the face integral, is
formally studied by Kronbichler et al.~\cite{kronbichler_fast_2019} in the
deal.II framework and many design choices are tested and compared. Recently this
was also studied in the Dune framework~\cite{bastian_matrix-free_2019} with some
different design choices, and excellent performance results were
reported~\cite{kempf_automatic_2020}. It is necessary to highlight the major
differences between our designs.

\begin{itemize}
    \item We target arbitrary unstructured meshes filled with various shape types,
    so we miss many optimisations specialised for hexahedra-only meshes, making it
    difficult to achieve the same performance as reported in other work.
    \item deal.II includes structured affine meshes consisting of identical
    hexahedral elements with shared Jacobian and geometric data. This has
    significantly lower memory loading requirements and increase operation intensity.
    \item deal.II also uses Hermite polynomials for the element basis, which
    allows the evaluation of derivatives on traces from two layers of points. On
    the contrary, for the element bases we listed in Table \ref{tab:elements &
    bases}, the whole element space is always needed to evaluate the derivative,
    which is theoretically more expensive.
    \item Dune uses a different vectorisation strategy: instead of evaluating 4
    elements at a time, they evaluate 1 solution and 3 derivatives at the same
    time, by the same sum-factorisation kernel. To fill the 512-bit register, they
    also need to pack data on different points, so the data layout is more
    complicated.
\end{itemize}

\subsection{Performance benchmarking}

The performance benchmarks are performed on two different CPU computing nodes. The older machine has two Xeon E5-2697 v4 Broadwell CPUs with AVX2, mostly for comparison with the previous studies. The newer machine has two EPYC 9554 CPUs with AVX512. The detailed specifications of the hardware and software environment are summarised in Table \ref{tab:hardware}. 

\begin{table}
    \caption{Hardware specifications and software configurations used for performance benchmarks}
    \centering
    \begin{tabular}{lll}
    \hline
        CPU Model  & Xeon E5-2697 v4 & EPYC 9554 \\
    \hline
        Architecture & Broadwell  & Zen4 \\
        SIMD capability  &  AVX2 (256bit) & AVX512 (512bit) \\
        Base clock speed & 2.3 GHz  & 3.1 GHz \\
        L2 cache per socket &  4.5 MB   & 64 MB  \\
        L3 cache per socket &  45 MB   & 256 MB  \\
        Memory           & 4-ch DDR4-2400 & 12-ch DDR5-4800 \\
        Cores per socket & 18      & 64   \\
        Socket per node  & 2       & 2     \\
        Peak GFLOPS per node & 1152  & 7611 \\
        Compiler and flags & GCC 12.2.0, -O3 & GCC 11.4.0, -O3 \\
        Parallel library &  Open MPI 4.1.4 & Open MPI 4.1.2 \\
        Profiler         & likwid 5.2.0 &  likwid 5.4.1 \\
    \hline
    \end{tabular}
    \label{tab:hardware}
\end{table}

\begin{figure}
    \centering
    \begin{subfigure}[b]{0.32\textwidth}
        \centering
        \includegraphics[width=\textwidth]{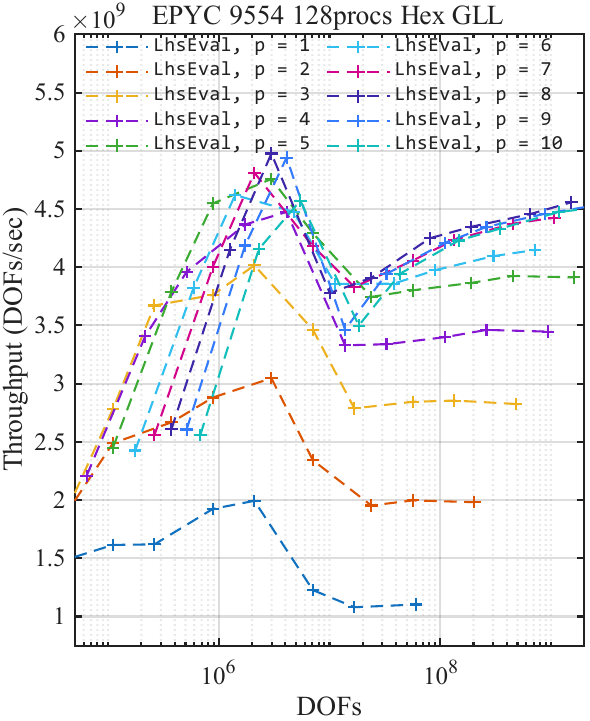}
        \caption{}
        \label{}
    \end{subfigure}
    \begin{subfigure}[b]{0.32\textwidth}
        \centering
        \includegraphics[width=\textwidth]{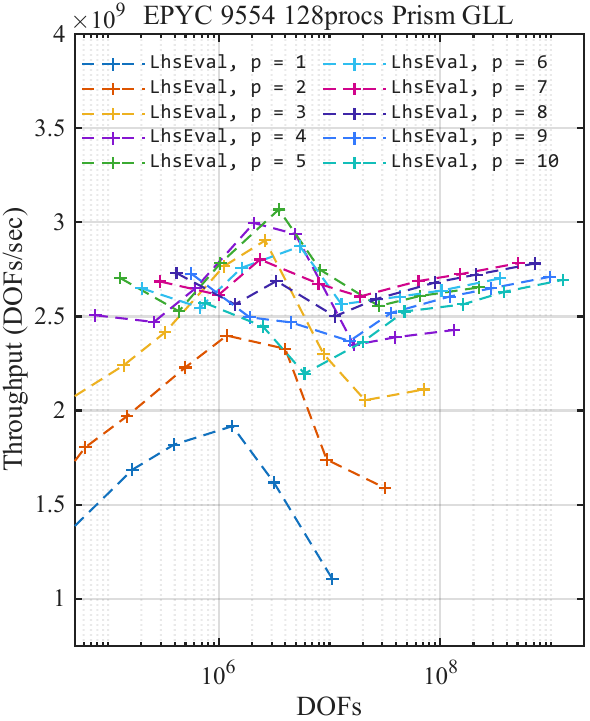}
        \caption{}
        \label{}
    \end{subfigure}
    \begin{subfigure}[b]{0.32\textwidth}
        \centering
        \includegraphics[width=\textwidth]{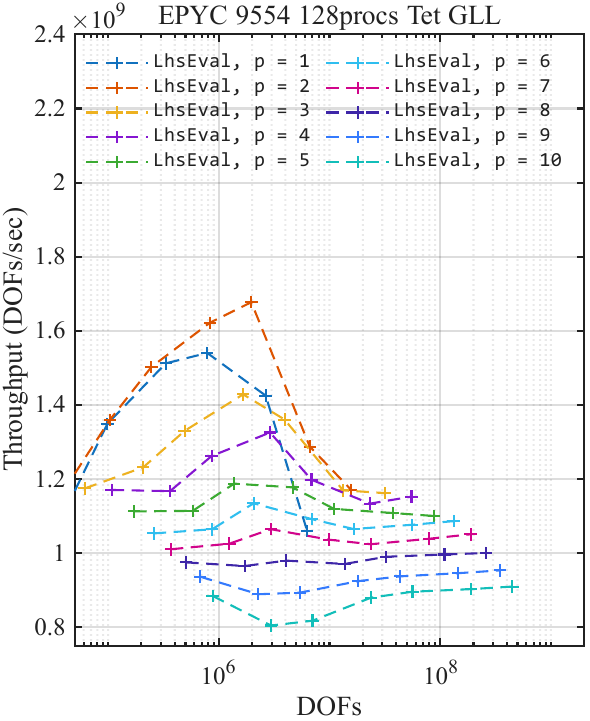}
        \caption{}
        \label{}
    \end{subfigure}
    \caption{The throughput versus problem sizes (DOFs) and element orders, tested on EPYC 9554. (a) hexahedrons; (b) Prisms; (c) Tetrahedrons.}
    \label{fig:TPUT-DOF}
\end{figure}

We focus on the maximum performance that can be achieved using all cores of a given node. Throughput is used to measure the performance, which is the number of degrees of freedom (DOF) processed per second by the target operator, calculated based on the elapsed time of the operator and the total DOFs it processes. Floating-point operations per second (FLOPS) can indicate how much potential we have exploited from the machine, which is obtained via the profiling tool suite \texttt{likwid}.

The throughput is affected by many factors: the spectral element shapes, basis
order and quadrature order directly determine the computing complexity, whilst
different problem sizes affect cache utilisation. Fig.~\ref{fig:TPUT-DOF} gives
an overview of the throughput of LHS evaluation for varying DOFs, element orders
and shapes. For brevity, we only test the modal basis functions with $p+2$ GLL
(or $p+1$ GR) points along each coordinate. All the meshes we used for
benchmarks are structured, such as $64^3$ and $48^3$, but the code treats them
as general unstructured meshes. No special techniques are applied during the
pre-processing to improve performance. As the problem size increases, the
throughput first increases to the peak at around $10^6$ DOFs, where the cache
performance is maximised; thereafter, the throughput experiences a short decline
and soon recovers and levels off. Cache blocking techniques ensure that the
memory access pattern is nearly independent of the data size, and so is the
performance. For hexahedral and prismatic elements, increasing $p$ from 1 to 5
significantly improves the performance, but for tetrahedral elements, increasing
the order reduces the performance. The reason is twofold. When we fix the
problem size and increase the order, the mesh will be sparser, but local element
operations are denser, leading to better performance. However, increasing the
order also increases the computing complexity, which should reduce the
performance. For tetrahedral elements, the second reason dominates, and their
operators cannot be optimised to the same level as hexahedral ones due to the
variable loop counts.

\begin{figure}
    \centering
    \begin{subfigure}[b]{0.32\textwidth}
        \centering
        \includegraphics[width=\textwidth]{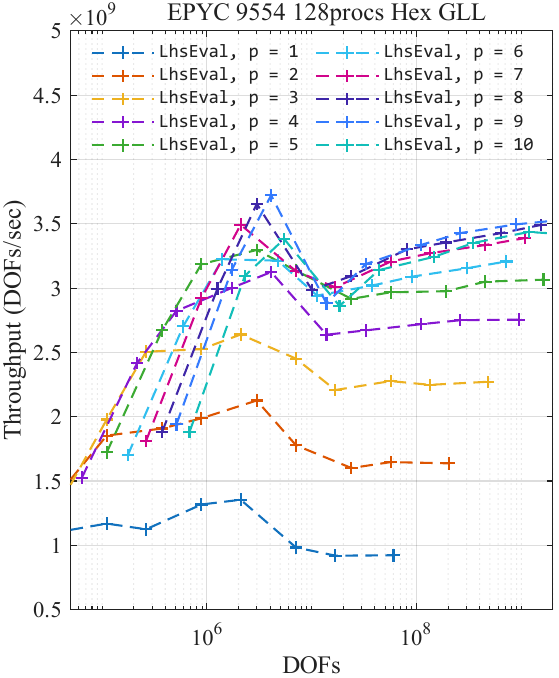}
        \caption{}
        \label{}
    \end{subfigure}
    \begin{subfigure}[b]{0.32\textwidth}
        \centering
        \includegraphics[width=\textwidth]{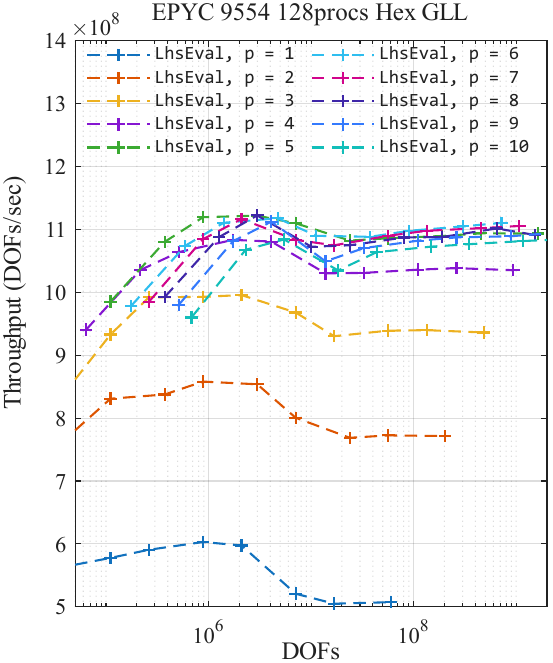}
        \caption{}
        \label{}
    \end{subfigure}
    \caption{Throughput vs. problem sizes (DOFs) and element order, tested on EPYC 9554, for simulated
      non-conformal interfaces via (a) interpolation or (b) \texttt{TraceIProduct} on element boundaries.}
    \label{fig:TPUT-DOF-forceinterp}
\end{figure}

A simple way to evaluate the performance in processing non-conforming interfaces
is to force the use of \texttt{TraceIProduct} or interpolation on all
interfaces, even if the discretisation is conforming. Fig.~\ref{fig:TPUT-DOF-forceinterp}
provides these results and can be directly compared with
Fig.~\ref{fig:TPUT-DOF}. We see that interpolation incurs around $30\%$
performance loss with respect to the conformal case, but it is still
significantly faster than \texttt{TraceIProduct}. One factor is that the existing
framework performs kernel fusing and cache blocking based on the element rather
than the trace. Thus, some optimisation cannot be applied to
\texttt{TraceIProduct} in the same way as \texttt{ScatrInterp}. In this work, we
primarily focus on the performance characteristics of interpolation approach.

\begin{figure}
    \centering
    \begin{subfigure}[b]{0.32\textwidth}
        \centering
        \includegraphics[width=\textwidth]{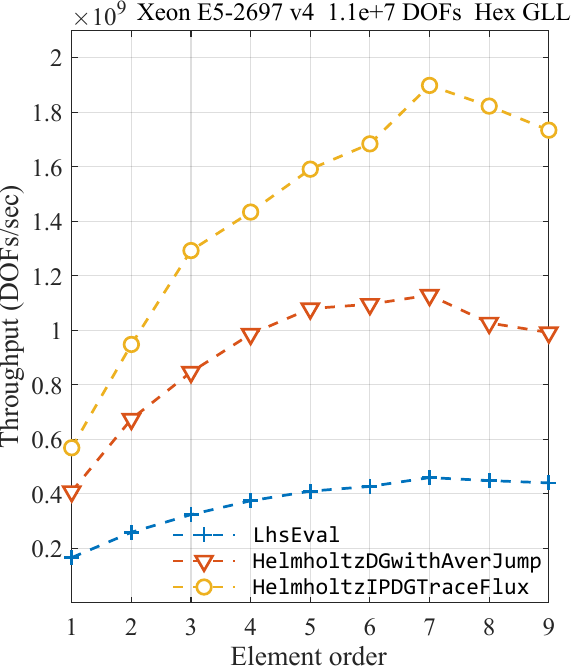}
        \caption{}
        \label{}
    \end{subfigure}
    \begin{subfigure}[b]{0.32\textwidth}
        \centering
        \includegraphics[width=\textwidth]{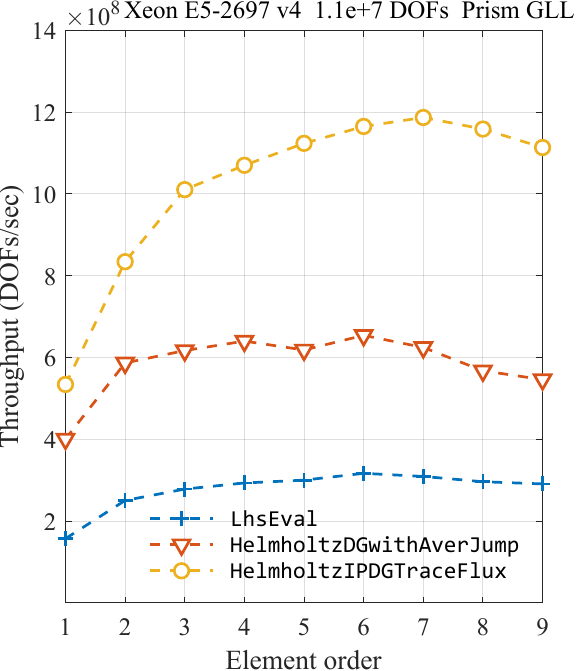}
        \caption{}
        \label{}
    \end{subfigure}
    \begin{subfigure}[b]{0.32\textwidth}
        \centering
        \includegraphics[width=\textwidth]{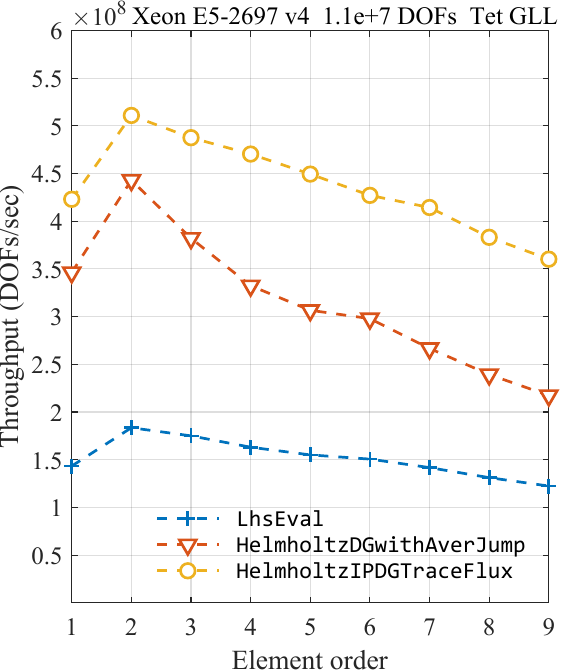}
        \caption{}
        \label{}
    \end{subfigure}
    \caption{The throughput versus element order, tested on Xeon E5-2697 v4, with fixed-size ($1.1\times10^{7}$ DOFs) hexahedral (a), prismatic (b) and tetrahedral (c) meshes.}
    \label{fig:TPUT-order}
\end{figure}

Most previous studies choose a single mesh for each order for the performance benchmark, which produces more concise results, highlighting the difference between orders. The problem size is significantly bigger than the L3 caches, reflecting typical workloads. For the Xeon node, the tested problem sizes are around $1\times10^{7}$ DOFs. The throughput results are shown in Fig. \ref{fig:TPUT-order}. The maximum throughput we achieved is about $4\times10^{8}$ DOFs/s on the Broadwell CPU, lower than the \textit{fusion} design given in \cite{kempf_automatic_2020}, but still comparable to their \textit{hybrid} design. The throughput of two sub-operators, \texttt{HelmholtzDGwithAverJump} and \texttt{HelmholtzIPDGTraceFlux}, gradually deviate from each other as the order increases. This is expected because the difference between element space and trace space grows larger for higher orders, and the volume flux also needs more operations to evaluate than the interface flux.

\begin{figure}
    \centering
    \begin{subfigure}[b]{0.33\textwidth}
        \centering
        \includegraphics[width=\textwidth]{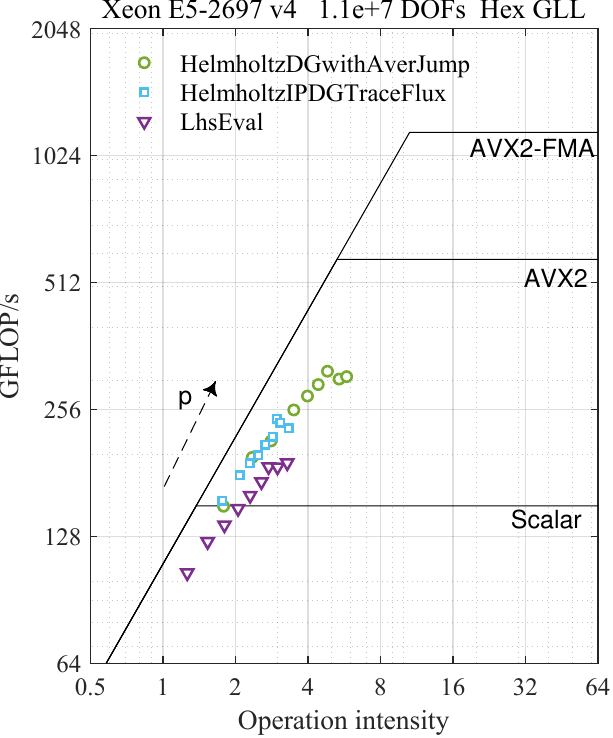}
        \caption{}
        \label{}
    \end{subfigure}
    \begin{subfigure}[b]{0.31\textwidth}
        \centering
        \includegraphics[width=\textwidth]{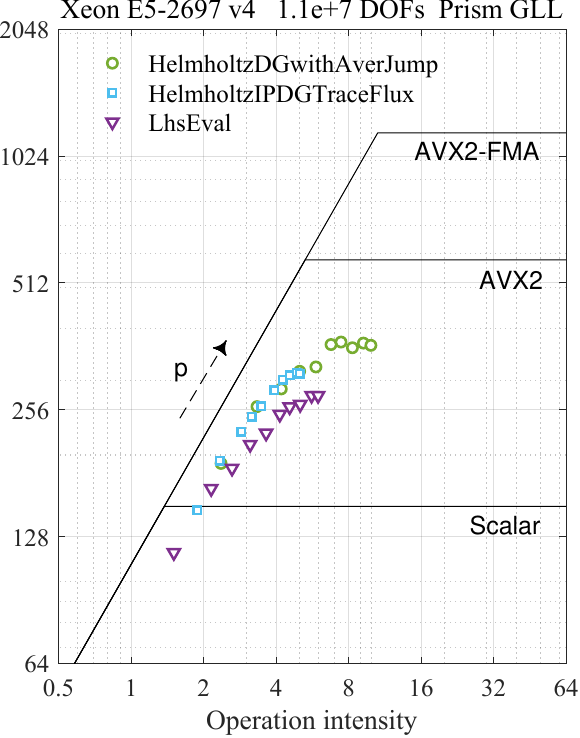}
        \caption{}
        \label{}
    \end{subfigure}
    \begin{subfigure}[b]{0.31\textwidth}
        \centering
        \includegraphics[width=\textwidth]{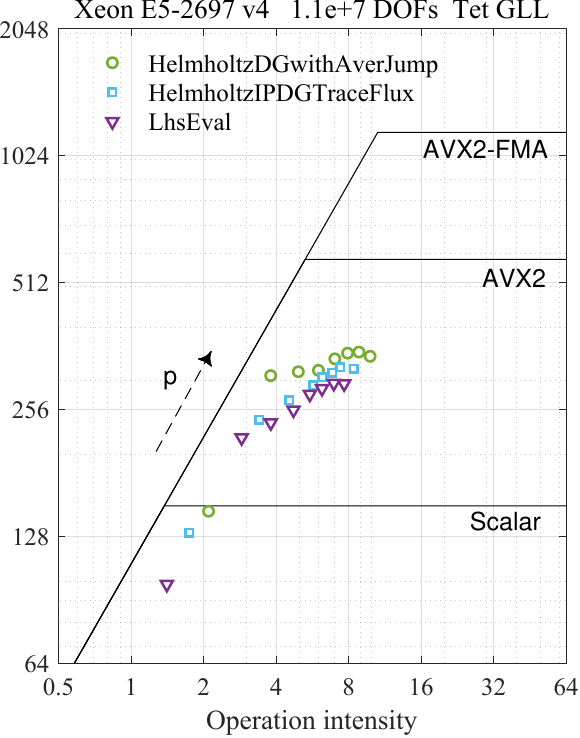}
        \caption{}
        \label{}
    \end{subfigure}
    \caption{The roofline results obtained on Xeon E5-2697 v4. The order increases from 1 to 10 as shown by the arrow. (a) Hexahedrons; (b) Prisms; (c) Tetrahedrons.}
    \label{fig:roofline-colo01}
\end{figure}

\begin{figure}
    \centering
    \begin{subfigure}[b]{0.33\textwidth}
        \centering
        \includegraphics[width=\textwidth]{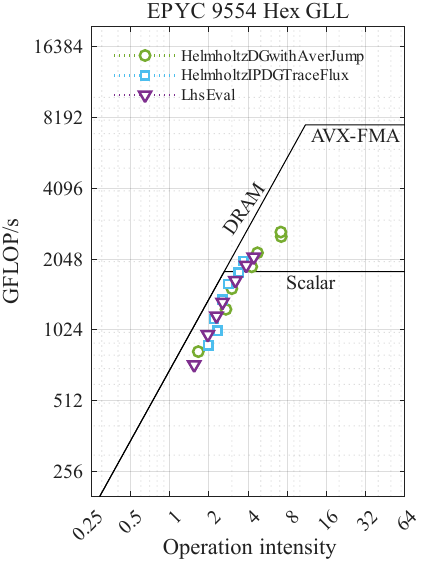}
        \caption{}
        \label{}
    \end{subfigure}
    \begin{subfigure}[b]{0.305\textwidth}
        \centering
        \includegraphics[width=\textwidth]{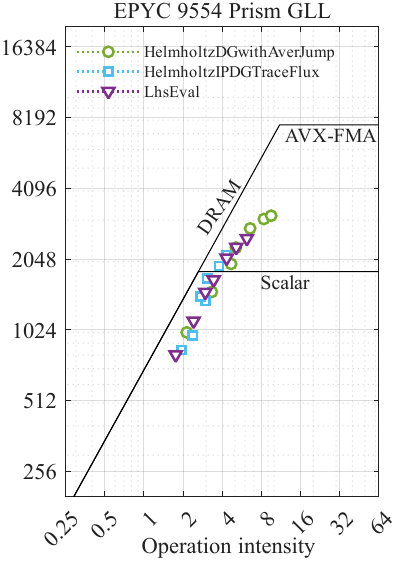}
        \caption{}
        \label{}
    \end{subfigure}
    \begin{subfigure}[b]{0.31\textwidth}
        \centering
        \includegraphics[width=\textwidth]{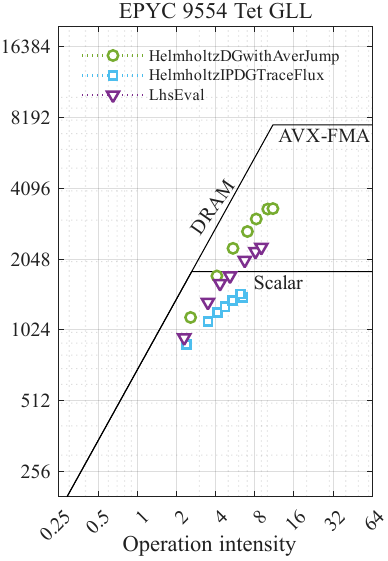}
        \caption{}
        \label{}
    \end{subfigure}
    \caption{The roofline results obtained on EPYC 9554. The orders are 1,2,3,4,6,8,10. (a) Hexahedrons; (b) Prisms; (c) Tetrahedrons.}
    \label{fig:roofline-kingscross}
\end{figure}

Finally, we present roofline results to show the utilisation of the
hardware. Our previous studies, e.g.~\cite{moxey_efficient_2020}, demonstrated
that the Helmholtz volume operations can reach AVX2 peak FLOPS on the same Xeon
node. As a comparison, we present roofline results on the same machine for the \texttt{LhsEval} and
two sub-operators in Fig.~\ref{fig:roofline-kingscross}.  The test case size is
consistently set to $1.1\times10^7$ DOFs for each order.  For this Xeon node, we
highlight 3 peak FLOPS rooflines obtained by \texttt{likwid-bench}.  The highest
corresponds to the theoretical peak of the CPU, which can only be achieved with
perfect use of AVX2 and fused-multiply add operations (FMA); the second highest
is the peak FLOPS achievable without FMA; and the lowest indicates the peak
level without the use of any vectorisation.  The two sub-operators fall between
the second and the third line, roughly 20-30\% of peak FLOPS, but still lower
than the results reported by~\cite{muthing_high-performance_2017}.  The roofline
results obtained on the EPYC node are displayed in
Fig.~\ref{fig:roofline-colo01}.  The difference between FMA, AVX512, and AVX2
peak FLOPS is marginal on this node owing to the implementation of AVX512 on
EPYC CPUs, resulting in only two distinct rooflines.  The tested problem size
exceeds $1\times10^8$ DOFs to exhaust memory and cache bandwidth.  Generally, as
the order increases, both FLOPS and operational intensity rise, yet they remain
evidently memory-bound.  Tetrahedral elements achieve slightly higher FLOPS and
also higher operation intensity than other shapes.  The operation intensity is
similar on the two machines, which is more related to the design itself.  We will
seek further optimisation in a later work, e.g., more aggresive kernel fusing
and tuning the cache block size, in order to increase the operation intensity.

\section{Conclusions}

In this work, we have discussed the evaluation of the interface flux terms for
DG methods, from formulation to implementation, aiming at a unified and low-cost
DG framework for general spectral elements of different choices of shape, basis
and quadrature. We have introduced an optimised matrix-free workflow for the DG
solver. The key idea of the workflow is that all terms should be evaluated in
the physical space, and only at the end do we transform physical space to
coefficient space. Unlike collocated nodal elements, the transformation between
coefficient space and physical space in a general spectral element is relatively
expensive. So in this workflow, we reduce the number of transform operations to
only two: a \texttt{BwdTrans} at the beginning, and an \texttt{IProduct} at the
end.

The key operators that bridge the traces and elements have been discussed in
detail. First, we identify two ways to get trace data from elements,
\texttt{TracePhysEval} and \texttt{GathrInterp}. They are initially designed for
Gauss and Gauss-Lobatto points, respectively, but can be generalised for all
cases if necessary, and they are numerically equivalent. We then identify two
ways to handle non-conforming interfaces, such as different polynomial orders
between adjacent elements. For SIPG methods, this is critical since the system
matrix can easily become non-symmetric if interface flux terms are inconsistent
across the non-conforming interfaces, which finally causes conjugate gradient
iterations to fail to converge. To resolve this issue, the first strategy is to
unify the quadrature and flux evaluation on both sides, by defining the shared
trace space. The second strategy is to use sufficient quadrature points so that
we can directly interpolate local element data to the adjacent side and evaluate
the flux locally, known as point-to-point interpolation.  The shared trace space
approach is mathematically equivalent to the mortar element method, but we only
copy data between physical spaces, instead of performing the $L^2$ projection
using the inverse mass matrix, making it much cheaper. In this work, we only
focus on the polynomial non-conforming instead of geometric non-conforming in
the discussion and numerical tests, although this extension is theoretically
possible for both approaches.

Finally, an initial matrix-free implementation based on the optimised workflow is
provided in detail. We exploit the performance on CPUs by improving cache data reuse and using explicit SIMD instructions. 
The interpolation approach via \texttt{ScatrInterp} has been proven to be much more efficient than 
performing \texttt{TraceIProduct} on shared traces, regardless of whether they are conforming or not.
This is an important resaon why we still emphasise the implementation of the interpolation approach in the current work, although it has numerical defects.
We aim at unstructured meshes with hybrid element types, so it is naturally more difficult to optimise. 
Although our design does not achieve the highest performance in comparison to other studies, 
we enabled an efficient design for various element types and also non-conforming meshes. 
We have provided performance benchmark results for various problem sizes, element orders, and shape types, which can be
a good reference for those interested in performance topics. Aside from optimising the matrix-free implementation, 
future work will also focus on large-scale parallel computing, with a particular focus on reducing communication overheads and effective preconditioning strategies.

\appendix \label{Apendix:A}
\section{Relationship between flux terms and symmetry of system matrix} 

Assume there are two overlapping elemental boundaries, $\Gamma_{r}$ and $\Gamma_{l}$, from two adjacent element $\Omega_{r}$ and $\Omega_{l}$. We use $r$ and $l$ to distinguish the right and left side objects. For instance, the bases of two elements are $\phi^r_{i}$, $\phi^l_{i}$ and coefficients are $\hat{u}^r_{i}$, $\hat{u}^l_{i}$. The unit normals always point outwards on the elemental boundary, so $\boldsymbol{n}^l = - \boldsymbol{n}^r$. With these notations, the three flux terms in Eq. \eqref{eq:IP} become
\begin{equation*}
    \text{Symmetry flux added to the left element : } -(\nabla \phi^l_{i} \cdot \boldsymbol{n}^{l}, \phi^l_{j})_{\Gamma_{l}} \frac{1}{2}\hat{u}^l_{j} + (\nabla \phi^l_{i} \cdot \boldsymbol{n}^{l}, \phi^r_{j})_{\Gamma_{l}} \frac{1}{2}\hat{u}^r_{j}
\end{equation*}
\begin{equation*}
    \text{Trace flux added to the left element : } -(\phi^l_{i}, \boldsymbol{n}^{l} \cdot \nabla \phi^l_{j})_{\Gamma_{l}} \frac{1}{2}\hat{u}^l_{j} - (\phi^l_{i}, \boldsymbol{n}^{l} \cdot \nabla \phi^r_{j})_{\Gamma_{l}} \frac{1}{2}\hat{u}^r_{j}
\end{equation*}
\begin{equation*}
    \text{Penalty flux added to the left element : } \tau(\phi^l_{i},  \phi^l_{j})_{\Gamma_{l}} \hat{u}^l_{j} - \tau(\phi^l_{i}, \phi^r_{j})_{\Gamma_{l}} \hat{u}^r_{j}
\end{equation*}

Flux added to the right side can be written down similarly. Here, we
deliberately perform the inner product in the local element boundaries to show
how it affects symmetry of the system. Now we collect contributions from the
left and right sides:

\begin{equation*} \label{eq:L2L}
\begin{aligned}
    & \quad\quad\quad\quad \mathbf{M}_A \quad \quad \quad \quad \quad \quad  \mathbf{M}'_A \quad\quad\quad\quad\quad\quad \mathbf{M}_E \\
    \text{Left to left :} &
    \left[ -\overbrace{ (\phi^l_{i}, \boldsymbol{n}^{l} \cdot \nabla \phi^l_{j})_{\Gamma_{l}} } - \overbrace{ (\nabla \phi^l_{i} \cdot \boldsymbol{n}^{l}, \phi^l_{j})_{\Gamma_{l}} } + \overbrace{ 2\tau(\phi^l_{i}, \phi^l_{j})_{\Gamma_{l}} } \right]\frac{1}{2}\hat{u}^l_{j} \\
    & \quad\quad\quad\quad \mathbf{M}_B \quad \quad \quad \quad \quad \quad  \mathbf{M}'_B \quad\quad\quad\quad\quad\quad \mathbf{M}_F \\
    \text{Right to right :} &
    \left[ -\overbrace{ (\phi^r_{i}, \boldsymbol{n}^{r} \cdot \nabla \phi^r_{j})_{\Gamma_{r}} } - \overbrace{ (\nabla \phi^r_{i} \cdot \boldsymbol{n}^{r} ,\phi^r_{j})_{\Gamma_{r}} } + \overbrace{ 2\tau(\phi^r_{i}, \phi^r_{j})_{\Gamma_{r}} } \right]\frac{1}{2}\hat{u}^r_{j} \\
\end{aligned}
\end{equation*}

\begin{equation*} \label{eq:L2R}
\begin{aligned}
    & \quad\quad\quad\quad \mathbf{M}_C \quad \quad \quad \quad \quad \quad  \mathbf{M}'_D \quad\quad\quad\quad\quad\quad \mathbf{M}_G \\
    \text{Right to left :} &
    \left[ -\overbrace{ (\phi^l_{i}, \boldsymbol{n}^{l} \cdot \nabla \phi^r_{j})_{\Gamma_{l}} } +\overbrace{ (\nabla \phi^l_{i} \cdot \boldsymbol{n}^{l}, \phi^r_{j})_{\Gamma_{l}} } - \overbrace{ 2\tau(\phi^l_{i}, \phi^r_{j})_{\Gamma_{l}} } \right]\frac{1}{2}\hat{u}^r_{j} \\
    & \quad\quad\quad\quad \mathbf{M}_D \quad \quad \quad \quad \quad \quad  \mathbf{M}'_C \quad\quad\quad\quad\quad\quad \mathbf{M}'_G \\
    \text{Left to right :} &
    \left[ -\overbrace{ (\phi^r_{i}, \boldsymbol{n}^r \cdot \nabla \phi^l_{j})_{\Gamma_{r}} } + \overbrace{ (\nabla \phi^r_{i} \cdot \boldsymbol{n}^r,\phi^l_{j})_{\Gamma_{r}} } -\overbrace{ 2\tau(\phi^r_{i}, \phi^l_{j})_{\Gamma_{r}} }\right]\frac{1}{2}\hat{u}^l_{j}
\end{aligned}
\end{equation*}

The locations of these elemental matrices in the system matrix are:
\begin{equation*}
    \begin{bmatrix} 
        -\mathbf{M}_A -  \mathbf{M}'_A + \mathbf{M}_E & \cdots & - \mathbf{M}_C + \mathbf{M}'_D - \mathbf{M}_G \\
                \cdots & \ddots & \cdots \\ 
        \mathbf{M}'_C - \mathbf{M}_D - \mathbf{M}'_G & \cdots & -\mathbf{M}_B - \mathbf{M}'_B + \mathbf{M}_F 
    \end{bmatrix}
\end{equation*}
Matrices A, B, E, and F are added to the diagonal entries of the blocked system matrix, while C, D, and G are added to the off-diagonal entries. To get a symmetric system, we must ensure 
\begin{equation*} \label{eq:local matrices}
    \mathbf{M}'_A = \mathbf{M}^T_A, \quad \mathbf{M}'_B = \mathbf{M}^T_B , \quad \mathbf{M}_E = \mathbf{M}^T_E , \quad \mathbf{M}_F = \mathbf{M}^T_F
\end{equation*}
and
\begin{equation*} \label{eq:inter-elemental matrices}
    \mathbf{M}'_C = -\mathbf{M}^T_C, \quad \mathbf{M}'_D = -\mathbf{M}^T_D  , \quad \mathbf{M}'_G = \mathbf{M}^T_G
\end{equation*}
The first relation is always satisfied. As for the second one, we must ensure that quadrature on the overlapped boundaries $\Gamma_l$ $\Gamma_r$ are either the same or have sufficient accuracy.

\bibliographystyle{unsrt}  
\bibliography{references}  

\end{document}